\documentclass[11pt, dvipsnames]{amsart}
\usepackage{amsthm,amssymb,amsmath}

\usepackage[left=3cm, right=3cm]{geometry}
\usepackage{microtype} 

\usepackage{mathtools} 
\usepackage[english]{babel}
\usepackage[shortlabels]{enumitem}
\usepackage{url}
\usepackage{hyperref}

\theoremstyle{plain}
\newtheorem{theorem}{Theorem}[section]
\newtheorem{lemma}[theorem]{Lemma}
\newtheorem{proposition}[theorem]{Proposition}
\newtheorem{corollary}[theorem]{Corollary}
\newtheorem{mainthm}{Main Theorem}

\theoremstyle{definition}
\newtheorem{definition}[theorem]{Definition}
\newtheorem{remark}[theorem]{Remark}
\newtheorem{example}[theorem]{Example}
\newtheorem{situation}[theorem]{Situation}

\def\C{\mathbb{C}}
\def\R{\mathbb{R}}
\def\Q{\mathbb{Q}}
\def\Z{\mathbb{Z}}
\def\F{\mathbb{F}}
\def\Fq{\mathbb{F}_q}

\DeclareMathOperator{\ord}{ord}
\DeclareMathOperator{\SIGN}{SIGN}
\DeclareMathOperator{\MULT}{MULT}

\newcommand{\Wq}[1]{\mathcal{W}_{q}({#1})}
\newcommand{\abs}[1]{\left\vert#1\right\vert}

\usepackage{color}

\title[Weil polynomials of small degree]{Weil polynomials of small degree}
\date{\today}

\author{Stefano Marseglia}
\address{Laboratoire Jean Alexandre Dieudonn\`e, Universit\`e C\^ote Azur, 06108 Nice Cedex 2, France}
\email{stefano.marseglia@univ-cotedazur.fr}

\keywords{
    abelian varieties over finite fields, roots of real polynomials
}

\subjclass{
    14K15,  
    12D10,  
    26C10. 
}

\begin{document}

\begin{abstract}
    Honda and Tate showed that the isogeny classes of abelian varieties of dimension $g$ over a finite field $\mathbb{F}_q$ are classified in terms of $q$-Weil polynomials of degree $2g$, that is, monic integer polynomials whose set of complex roots consists of $g$ conjugate pairs of absolute value $\sqrt{q}$.
    There are descriptions of the space of such polynomials for $g \leq 5$, but for $g=3$, $4$ and $5$, these results contain mistakes.
    We correct these statements.
    Our proofs build on a criterion that determines when a real polynomial has only real roots in terms of the non-necessarily distinct roots of its first derivative.
\end{abstract}

\maketitle

\section{Introduction}\
\label{sec:intro}
\subsection*{Abelian varieties over finite fields up to isogeny.}
Let $A$ be an abelian variety of dimension $g$ over a finite field $\Fq$ of positive characteristic $p$.
Let $h_A(x)$ be the characteristic polynomial of the $\F_q$-Frobenius endomorphism of $A$ acting of the  $\ell$-adic Tate module $T_\ell(A)$ of~$A$, where $\ell$ is any rational prime coprime to $q$.
The polynomial $h_A(x)$ does not depend on the choice of $\ell$ and is a \emph{$q$-Weil polynomial}, that is, a monic polynomial with integer coefficients, of even degree $2g$, and whose set of complex roots has the form $\{w_1,\bar{w_1},\ldots,w_g,\bar{w_g} \}$, where each~$w_i$ has complex absolute value $\sqrt{q}$; see Definition~\ref{def:qWeil}.

Remarkably, the characteristic polynomial $h_A(x)$ completely determines the $\Fq$-isogeny class of $A$.
Indeed, Tate proved that two abelian varieties $A$ and $B$ are $\Fq$-isogenous if and only if $h_A(x)=h_B(x)$; see \cite[Theorem 1.(c)]{Tate66}.
In other words, the association $A\mapsto h_A(x)$ induces an injective map from the set of $\F_q$-isogeny classes of abelian varieties defined over $\F_q$ to the set of $q$-Weil polynomials.
This map is not quite surjective, but celebrated results of Honda and Tate describe its image; see \cite{Tate66,Honda68,Tate69} or Remark~\ref{rmk:HTtheory} for details.

\subsection*{Computing \texorpdfstring{$q$}{q}-Weil polynomials.}
One upshot of the discussion above is that to classify and enumerate the $\Fq$-isogeny classes of abelian varieties over $\Fq$ of dimension $g$, it is desirable to describe or compute the set $\Wq{g}$ of $q$-Weil polynomials of degree $2g$.

An efficient algorithm to compute $\Wq{g}$, for any $g$, is available at \cite{Kedlaya_repo}.
It has been included in SageMath \cite{sagemath} since January 2020.
The algorithm originated in the work of Abbott, Kedlaya and Roe \cite{AbbottKedlayaRoe10} and was further developed by Kedlaya in \cite{Kedlaya08} and Kedlaya-Sutherland in \cite{KedlayaSutherland16}.
It was then used by Dupuy, Kedlaya, Roe and Vincent to tabulate and study the isogeny classes of abelian varieties of dimension $g$ over $\Fq$ for various values of $g$ and $q$; see \cite{LMFDB_art_21}.
The result of their computations is available on the LMFDB \cite{lmfdb}.

Next to such computational results, there has been work to characterize the coefficients of the polynomials in $\Wq{g}$ for small values of $g$.
The task is trivial for $g=1$; see Example~\ref{ex:g1}.
Such a description is also well-known for $g=2$; see for example Maisner and Nart \cite[Lemma~2.1]{MaisnerNart02}, where it is deduced from work of R\"uck \cite[Lemma~3.1]{Ruck90}, or Proposition~\ref{ex:g2} below.

Analogous descriptions of the sets $\Wq{3}$, $\Wq{4}$ and $\Wq{5}$ were provided by Haloui in \cite{Haloui10}, Haloui and Singh in \cite{HalouiSingh12}, and Sohn in \cite{Sohn13}, respectively, but, as pointed out in \cite{LMFDB_art_21}, they all contain mistakes of various severity.
In Bradford's thesis \cite{Bradford12}, there is an attempt to correct the statement for $\Wq{4}$, but it still contains a small mistake.
The main purpose of this paper is to correct all statements.
This is done in Main Theorems~\ref{mainthm:dim3}, \ref{mainthm:dim4} and \ref{mainthm:dim5}.
See Remarks~\ref{rmk:mainthmB}, \ref{rmk:mainthmC} and \ref{rmk:mainthmD} for the comparisons with the original results.

\subsection*{Our contributions and structure of the paper.}
The main idea to compute $\Wq{3}$, $\Wq{4}$ and $\Wq{5}$, which is shared by \cite{Haloui10}, \cite{HalouiSingh12} and \cite{Sohn13}, is to reduce the question of determining whether a monic integral polynomial $h(x)$ of degree $2g$ is a $q$-Weil polynomial to the determination of 
whether two sets $\mathcal{S^+}$ and $\mathcal{S^-}$ of $g$ complex numbers constructed from the $2g$ complex roots of $h(x)$ are subsets of  $\R_{\geq 0}$.
Section~\ref{sec:qWeil} is devoted to introduce $q$-Weil polynomials and prove this result; see Proposition~\ref{prop:plus_minus}.

In order to determine whether the roots of a real polynomial $f(x)$ of degree $g$ are real, one can study the roots of its first derivative $f'(x)$; this is done in Section~\ref{sec:poly_real_pos_roots}, in particular see Proposition~\ref{prop:real_pos_only}.
Even if the proposition is probably well-known, we could not locate it in the literature.
We provide a short proof in Section~\ref{sec:poly_real_pos_roots} and a longer but more direct proof in Appendix~\ref{sec:appendix}.

Taking derivatives lowers the degree, and there are formulas in terms of radicals for the roots of polynomials of degree $\leq 4$.
Therefore, if our polynomial $f(x)$ has degree $\leq 5$, then a recursive application of Proposition~\ref{prop:real_pos_only} allows us to determine whether all the roots of $f(x)$ are real solely in terms of its coefficients.
This idea is developed in Section~\ref{sec:poly_low_degree} leading to Propositions~\ref{prop:deg3} and \ref{prop:deg4} and Corollary~\ref{cor:monic_deg5}.

Finally, we combine these results with Proposition~\ref{prop:plus_minus} to prove Main Theorems~\ref{mainthm:dim3}, \ref{mainthm:dim4} and \ref{mainthm:dim5} stated below, whose proofs can be found in Sections~\ref{sec:deg6}, \ref{sec:deg8} and \ref{sec:deg10}, respectively.

In the rest of the paper, we will use the following conventions:
By positive and non-negative we mean $>0$ and $\geq 0$, respectively.
If $r$ is a positive real number then we denote the positive real $n$-th root of $r$ by $\sqrt[n]{r}$ or $r^{1/n}$ ; in all other cases, the same symbols denote an arbitrary $n$-th complex root.

\begin{mainthm} \label{mainthm:dim3}
    Let $h(x) \coloneqq x^6 + a_1 x^5 + a_2 x^4 + a_3 x^3 + a_2 q x^2 + a_1 q^2 x + q^3$ be a polynomial with integer coefficients. 
    The polynomial $h(x)$ belongs to $\Wq{3}$ if and only if all the following conditions hold:
    \begin{enumerate}[(a)]
        \item \label{mainthm:dim3:a} $|a_1| \leq 6 \sqrt{q}$,
        \item \label{mainthm:dim3:b_star} $4 \sqrt{q} |a_1 | - 9q \leq a_2 $
        \item \label{mainthm:dim3:b} $ a_2 \leq \dfrac{1}{3}a_1^2 + 3q$,
        \item \label{mainthm:dim3:c} $- \dfrac{2 }{27}a_1^3 + \dfrac{1}{3}a_1 a_2 + q a_1 - \dfrac{2}{27} (a_1^2 - 3 a_2 + 9q)^{3/2} \leq a_3 \leq - \dfrac{2}{27}a_1^3 + \dfrac{1}{3}a_1 a_2 + q a_1 + \dfrac{2}{27} (a_1^2 - 3 a_2 + 9q)^{3/2}$,
        \item \label{mainthm:dim3:d} $-2 q a_1 - 2 \sqrt{q} a_2 - 2 q \sqrt{q} \leq a_3 \leq -2 q a_1 + 2 \sqrt{q} a_2 + 2 q \sqrt{q}$,
    \end{enumerate}
    where the quantity between parentheses in \ref{mainthm:dim3:c} is a non-negative real number by \ref{mainthm:dim3:b}.
    Moreover, $h(x) \in \Wq{3}$ has a real root if and only if any of the inequalities 
    in \ref{mainthm:dim3:a}, \ref{mainthm:dim3:b_star} or \ref{mainthm:dim3:d}
    is an equality.
    Explicitly, $h(x)$ is in $\Wq{3}$ and has a real root if and only if there are non-negative integers $k$ and $\ell$ such that one of the following mutually exclusive cases occurs:
    \begin{enumerate}[(I)]
        \item \label{mainthm:dim3:rr0} $q$ is not a square and $h(x) = (x^2 - q)^2 h_0(x)$ for $h_0(x)\in \Wq{1}$.
        \item \label{mainthm:dim3:rr1} $q$ is a square, $k+\ell = 3$ and $h(x) = (x + \sqrt{q})^{2k} (x - \sqrt{q})^{2 \ell}$.
        \item \label{mainthm:dim3:rr2} $q$ is a square, $k+\ell = 2$ and $h(x) = (x + \sqrt{q})^{2k} (x - \sqrt{q})^{2 \ell} h_0(x)$ for $h_0(x) \in \Wq{1}$ with no real roots.
        \item \label{mainthm:dim3:rr3} $q$ is a square, $k+\ell = 1$ and $h(x) = (x + \sqrt{q})^{2k} (x - \sqrt{q})^{2 \ell} h_0(x)$ for $h_0(x)\in \Wq{2}$ with no real roots.
    \end{enumerate}
\end{mainthm}

\begin{remark}\label{rmk:mainthmB}
    Main Theorem~\ref{mainthm:dim3} is a revised version of \cite[Theorem~1.1]{Haloui10}.
    By comparing our results with the original one, we see that the author of \cite{Haloui10} has overlooked some cases with real roots, all occurring when $q$ is a square.
\end{remark}

\begin{mainthm}\label{mainthm:dim4}
    Let $h(x) \coloneqq x^8 + a_1 x^7 + a_2 x^6 + a_3 x^5 + a_4 x^4 + a_3 q x^3 + a_2 q^2 x^2 + a_1 q^3 x + q^4$ be a polynomial with integer coefficients.
    Consider the following inequalities:
    \begin{enumerate}[(a)]
        \item \label{mainthm:dim4:a} $\abs{a_1} \leq 8 \sqrt{q}$,
        \item \label{mainthm:dim4:b_star} $6 \sqrt{q} \abs{a_1} - 20 q \leq a_2$,
        \item \label{mainthm:dim4:b} $a_2 \leq \dfrac{3}{8}a_1^2 + 4q$,
        \item \label{mainthm:dim4:c} $- 9q a_1 -4 \sqrt{q} a_2 - 16 q\sqrt{q} \leq a_3 \leq - 9 q a_1 + 4 \sqrt{q} a_2  + 16 q\sqrt{q}$,
        \item \label{mainthm:dim4:d} $\dfrac{1}{2}a_1 a_2 - \dfrac{1}{8}a_1^3 + q a_1 - \dfrac{1}{216} (9 a_1^2 - 24 a_2 + 96 q)^{3/2} \leq a_3 \leq \dfrac{1}{2}a_1 a_2 - \dfrac{1}{8}a_1^3 + q a_1 + \dfrac{1}{216} (9 a_1^2 - 24 a_2 + 96 q)^{3/2}$,
        \item \label{mainthm:dim4:e} $2 \sqrt{q}|q a_1 + a_3| - 2 q a_2 - 2q^2 \leq a_4$,
    \end{enumerate}
    where the quantity between parentheses in \ref{mainthm:dim4:d} is a non-negative real number by \ref{mainthm:dim4:b}.
    Define the real numbers
    \[ 
        u_2 \coloneqq -\frac{3 a_1^2}{16} + \frac{a_2}{2} - 2q, \qquad
        u_3 \coloneqq -\frac{a_1^3}{32} + \frac{a_1 a_2}{8} + \frac{a_1 q}{4} - \frac{a_3}{4} \qquad\text{and}\qquad
        \Delta \coloneqq u_3^2 + \frac{4u_2^3}{27},
    \]
    and the set
    \[ S \coloneqq \left\{-\zeta^{k} \omega \left( \frac{9 u_3}{2} + \frac{3\sqrt{\Delta}}{2}\right) - \zeta^{-k} \overline{\omega} \left( \frac{9 u_3}{2} - \frac{3\sqrt{\Delta}}{2} \right)  + \frac{2 u_2^2}{3} : 0\leq k \leq 2 \right\}, \]
    where $\zeta$ is a primitive third root of unity and $\omega$ is a third root of $ \frac{1}{2}(- u_3 + \sqrt{\Delta})$.
    If $h(x) \in \Wq{4}$ or \ref{mainthm:dim4:a}--\ref{mainthm:dim4:d} hold then the set $S$ consists of real numbers.
    Assume this is the case, sort the elements of $S$ as $\theta_1 \leq \theta_2 \leq \theta_3$ and consider the inequality:
    \begin{enumerate}[(a)]\setcounter{enumi}{6}
        \item \label{mainthm:dim4:f} 
        $\dfrac{3}{256}a_1^4 - \dfrac{1}{16}a_1^2 a_2 - \dfrac{1}{2}qa_1^2 + \dfrac{1}{4}a_1 a_3 + 2 q a_2 - 2q^2 + \theta_1 \leq a_4 \leq \dfrac{3}{256} a_1^4 - \dfrac{1}{16}a_1^2 a_2 - \dfrac{1}{2}qa_1^2 + \dfrac{1}{4}a_1 a_3 + 2 q a_2 - 2q^2 + \theta_2$.
    \end{enumerate}
    Then:
    \begin{center}
        $h(x)$ belongs to $\Wq{4}$ if and only if \ref{mainthm:dim4:a}--\ref{mainthm:dim4:f} hold.
    \end{center}
    Moreover, $h(x) \in \Wq{4}$ has a real root if and only if any of the inequalities 
    in \ref{mainthm:dim4:a}, \ref{mainthm:dim4:b_star}, \ref{mainthm:dim4:c} or \ref{mainthm:dim4:e}
    is an equality.
    Explicitly, $h(x)$ is in $\Wq{4}$ and has a real root if and only if one of the following mutually exclusive cases occurs:
    \begin{enumerate}[(I)]
        \item \label{mainthm:dim4:rr0} $q$ is a square and, $h(x) = (x + \sqrt{q})^{2} h_0(x)$ or $h(x) = (x - \sqrt{q})^{2} h_0(x)$ with $h_0(x)\in \Wq{3}$.
        \item \label{mainthm:dim4:rr1} $q$ is not a square and $h(x) = (x^2 - q)^2 h_0(x)$ with $h_0(x) \in \Wq{2}$.
    \end{enumerate}
\end{mainthm}
\begin{remark}\label{rmk:mainthmC}
    An incorrect version of Main Theorem~\ref{mainthm:dim4} is given in \cite[Theorem~1.1]{HalouiSingh12}.
    The main issue is that the authors seem to have overlooked that one needs to carefully sort the elements of the set $S$ to define the bounds on $a_4$.
    In \cite{Bradford12}, the author provides a correction, with no proof, but it still contains a small typo in the definition of the constant $\omega$.
    The method of the proof that we present in Section~\ref{sec:deg8} is similar in spirit to one used in \cite{Haloui10} for $\Wq{3}$ and slightly differs from the one given in \cite{HalouiSingh12}.
\end{remark}

\begin{mainthm}\label{mainthm:dim5}
    Let 
    \[ h(x) \coloneqq x^{10} + a_1 x^9 + a_2 x^8 + a_3 x^7 + a_4 x^6 + a_5 x^5 + a_4 q x^4 + a_3 q^2 x^3 + a_2 q^3 x^2 + a_1 q^4 x + q^5 \]
    be a polynomial with integer coefficients.
    Consider the following inequalities, where the quantity between parentheses in \ref{eq:mainthm:dim5:c} is a non-negative real number by \ref{eq:mainthm:dim5:b}, 
    and the definitions of $\theta_1$ and $\theta_2$ in \ref{eq:mainthm:dim5:e} and of $\Lambda_1$ and $\Lambda_2$ in \ref{eq:mainthm:dim5:h} are given below:
    \begin{enumerate}[(a)]
        \item  \label{eq:mainthm:dim5:a}$|a_1| \leq 10 \sqrt{q}$
        \item  \label{eq:mainthm:dim5:b_star} $8 \sqrt{q} |a_1| - 35 q \leq a_2$,
        \item  \label{eq:mainthm:dim5:b} $a_2 \leq \dfrac{2}{5} a_1^2 + 5 q$,
        \item  \label{eq:mainthm:dim5:c} $- \dfrac{4}{25} a_1^3 + \dfrac{3}{5} a_1 a_2 + q a_1 - \dfrac{1}{50} (4 a_1^2 + 50q - 10 a_2)^{3/2} \leq a_3 \leq - \dfrac{4}{25} a_1^3 + \dfrac{3}{5} a_1 a_2 + q a_1 + \dfrac{1}{50} (4 a_1^2 + 50q - 10 a_2)^{3/2},$
        \item  \label{eq:mainthm:dim5:d}$-6 \sqrt{q} a_2 - 20 q a_1 - 50 q \sqrt{q} \leq a_3 \leq 6 \sqrt{q} a_2 - 20 q a_1 + 50 q \sqrt{q},$
        \item \label{eq:mainthm:dim5:e} $\dfrac{3 a_1^4}{125} - \dfrac{3 a_1^2 a_2}{25} - q a_1^2 + \dfrac{2 a_1 a_3}{5} + 3 q a_2 - 5 q^2 + 5\theta_1 \leq a_4 \leq \dfrac{3 a_1^4}{125} - \dfrac{3 a_1^2 a_2}{25} - q a_1^2 + \dfrac{2 a_1 a_3}{5} + 3 q a_2 - 5 q^2 + 5\theta_2,$
        \item  \label{eq:mainthm:dim5:f} $ 4 \sqrt{q}|4q a_1 + a_3| - 9 q a_2 - 25 q^2 \leq a_4,$
        \item  \label{eq:mainthm:dim5:g}  
        $-2 \sqrt{q} a_4 - 2q a_3 - 2q \sqrt{q} a_2 - 2q^2 a_1 - 2q^2 \sqrt{q} \leq a_5 \leq 2 \sqrt{q} a_4 - 2q a_3 + 2 q \sqrt{q} a_2 - 2q^2 a_1 + 2q^2 \sqrt{q}$.
        \item  \label{eq:mainthm:dim5:h} $\Lambda_1 + A \leq a_5 \leq \Lambda_2 + A$,
    \end{enumerate}
    with
    \[ A\coloneqq - \frac{4a_1^5}{3125} + \frac{a_1^3\left(a_2 + 15q \right)}{125} - \frac{a_1^2a_3}{25} -\frac{a_1\left( 3a_2q -a_4 +5q^2\right)}{5} + 2 q a_3.  \]
    Define
    \begin{align*}
        u_2 &\coloneqq - \frac{3 a_1^2}{25} + \frac{3 a_2}{10} - \frac{3q}{2}, \\
        u_3 &\coloneqq \frac{2 a_1^3}{125} - \frac{3 a_1 a_2}{50} - \frac{q a_1}{10} + \frac{a_3}{10}, \\
        u_4 &\coloneqq - \frac{3 a_1^4}{625} + \frac{3 a_1^2 a_2}{125} + \frac{q a_1^2}{5} - \frac{2 a_1 a_3}{25} - \frac{3 q a_2}{5} + \frac{a_4}{5} + q^2,\\
        \Delta &\coloneqq u_3^2 + \frac{4}{27} u_2^3,  \text{ and}\\
        S & \coloneqq \left\{-\zeta^{k} \omega \left( \frac{9 u_3}{2} + \frac{3\sqrt{\Delta}}{2}\right) - \zeta^{-k} \overline{\omega} \left( \frac{9 u_3}{2} - \frac{3\sqrt{\Delta}}{2} \right)  + \frac{2 u_2^2}{3} : 0\leq k \leq 2 \right\},
    \end{align*}
    where $\zeta$ is a primitive third root of unity and $\omega$ is a third root of $\frac{1}{2} (- u_3 + \sqrt{\Delta})$.
    If $h(x)$ is in $\Wq{5}$ or \ref{eq:mainthm:dim5:a}--\ref{eq:mainthm:dim5:d} hold then the elements in $S$ are real numbers $\theta_1 \leq \theta_2 \leq \theta_3$.
    If $u_3 = 0$, let $x_{i_1, i_2}$ for $i_1, i_2 \in \{+1, -1 \}$ be
    \begin{equation*}
        x_{i_1, i_2} \coloneqq i_1 \sqrt{- u_2 + i_2 \sqrt{u_2^2 - u_4}}.
    \end{equation*}
    If $u_3 \neq 0$, let $x_{i_1, i_2}$ for $i_1, i_2 \in \{+1, -1 \}$ be
    \begin{equation*}
        x_{i_1, i_2} \coloneqq \frac{i_1 \sqrt{2y} + i_2 \sqrt{- 4 u_2 - 2y - i_1 \frac{8 u_3}{\sqrt{2y}}}}{2},
    \end{equation*}
    where:
    \begin{align*}
        & v_2 \coloneqq - \frac{u_2^2}{3} - u_4, \\
        & v_3 \coloneqq \frac{2 u_2 u_4}{3} - \frac{2 u_2^3}{27} - 2 u_3^2, \\
        & \text{$C$ is a third root of }
        \begin{cases}
            -v_3 & \text{if $v_2=0$},\\
            \frac{1}{2}\left(- v_3 + \sqrt{v_3^2 + \frac{4}{27} v_2^3}\right) & \text{if $v_2\neq 0$},
        \end{cases}\\
        & y \coloneqq C - \frac{v_2}{3C} - \frac{2 u_2}{3}.
    \end{align*}

    If $h(x)$ is in $\Wq{5}$ or \ref{eq:mainthm:dim5:a}--\ref{eq:mainthm:dim5:f} hold then the four elements $x_{i_1, i_2}$ are real numbers $\gamma_1 \leq \gamma_2 \leq \gamma_3 \leq \gamma_4$, which we use to define 
    \[
        \Lambda_1 \coloneqq \max\left\lbrace H(\gamma_1), H(\gamma_3)\right\rbrace
        \qquad \text{and} \qquad
        \Lambda_2 \coloneqq \min \left\lbrace H(\gamma_2), H(\gamma_4)\right\rbrace,
    \]
    where 
    \[ H(x)\coloneqq -x^5 -\frac{10}{3}u_2 x^3 -10u_3x^2 - 5u_4x. \]
    Then:
    \begin{center}
        $h(x)$ is in $\Wq{5}$ if and only if  \ref{eq:mainthm:dim5:a}--\ref{eq:mainthm:dim5:h} hold.
    \end{center}
    Moreover, $h(x) \in \Wq{5}$ has a real root if and only if any of the inequalities 
    in \ref{eq:mainthm:dim5:a}, \ref{eq:mainthm:dim5:b_star}, \ref{eq:mainthm:dim5:d}, \ref{eq:mainthm:dim5:f} or \ref{eq:mainthm:dim5:g}
    is an equality.
    Explicitly, $h(x)$ is in $\Wq{5}$ and has a real root if and only if one of the following mutually exclusive cases occurs:
    \begin{enumerate}[(I)]
        \item \label{mainthm:dim5:rr0} $q$ is a square and, $h(x) = (x + \sqrt{q})^{2} h_0(x)$ or $h(x) = (x - \sqrt{q})^{2} h_0(x)$ with $h_0(x)\in \Wq{4}$.
        \item \label{mainthm:dim5:rr1} $q$ is not a square and $h(x) = (x^2 - q)^2 h_0(x)$ with $h_0(x) \in \Wq{3}$.
    \end{enumerate}
\end{mainthm}
\begin{remark}\label{rmk:mainthmD}
    Main Theorem~\ref{mainthm:dim5} should be compared with \cite[Theorem 2.1]{Sohn13}.
    One notable difference is in the bound for $a_3$ where the author is missing a term $50 q$.
    Moreover, it seems that the polynomials in $\Wq{5}$ with real roots are excluded from the description.
    Similar to \cite{HalouiSingh12}, the author seems to have overlooked that one needs to carefully sort the elements of the set $S$ to define the bounds for $a_4$.
\end{remark}

\begin{remark}\label{rmk:compare_w_kedlaya}
    We have turned our results into algorithms to compute $\Wq{g}$ for $g=1,\ldots,5$ and implemented them in the computer algebra system Magma \cite{magma}.
    The output of our implementation has been tested against the one produced using \cite{Kedlaya_repo} for a variety of prime powers $q$.
    We have found no discrepancies.
    Our implementation and the code to reproduce the tests can be found at \cite{test_code_repo}.
    We stress that our code is for testing purposes since it is slower than \cite{Kedlaya_repo}.
    The reason is that sometimes in one of the nested loops, say after we have candidates for $a_1,\ldots,a_i$, the interval that should contain $a_{i+1}$ turns out to be empty.
    The approach developed in \cite{Kedlaya_repo} specifically targets this issue by testing more conditions on $a_1,\ldots,a_i$ to minimize the chances of such a dead-end.
    See \cite[Section~3.1]{LMFDB_art_21} for details on this matter.
\end{remark}

\subsection*{Additional related literature}
We conclude the introduction by mentioning further work related to $q$-Weil polynomials.
As anticipated above and explained in detail in Remark~\ref{rmk:HTtheory}, not every $q$-Weil polynomial of degree $2g$ is the characteristic polynomial of Frobenius of some abelian variety of dimension $g$.
Explicit conditions for this to happen are known for small values of $g$.
For $g=1$ this can be found in Waterhouse \cite[Section~4.1]{Waterhouse69}.
For $g=2$, we refer to R\"uck \cite[Theorem~1.1]{Ruck90}, Xing \cite[Lemma~1]{Xing94} and \cite[Lemma~1]{Xing96}, Maisner-Nart~\cite[Theorem~2.9]{MaisnerNart02}, or Howe-Nart-Ritzenthaler \cite[Section~14]{HoweNartRitzenthaler09}.
For $g=3 $ and $4$, see Xing \cite[Theorems~1 and 2]{Xing94}.
The characteristic polynomials of abelian varieties of dimension $g=5$ are singled out from $\Wq{5}$ in Hayashida \cite[Theorem~1.3]{Hayashida19}.
The same is done for $g=6$ in Sohn~\cite[Theorem~1.2]{Sohn20}, but note that in this case we don't have an explicit description of $\Wq{6}$.
Further partial results for arbitrary dimensions can be found in Hayashida \cite{Hayashida21}.

Isogeny classes of abelian varieties can be grouped according to their $p$-rank or Newton polygons.
It is well known that only two Newton polygons can occur for dimension $1$: the ordinary one ($p$-rank $1$) and the supersingular one ($p$-rank $0$).
In higher dimension the situation is more interesting. 
Classifications can be found:
in R\"uck \cite[Proof of Lemma~3.2]{Ruck90} for $g=2$; 
in Haloui \cite[Section~4]{Haloui10} and Xing \cite[Proposition~4]{Xing94} for $g=3$;
in Haloui-Singh \cite[Section~3]{HalouiSingh12} for $g=4$ and irreducible characteristic polynomials;
in \cite[Proof of Theorem~1.3]{Hayashida19} for $g=5$; in Sohn \cite[Section~2]{Sohn20} for $g=6$.

Another line of research is studying which characteristic polynomials give rise to an isogeny class that contains a principally polarized abelian variety or the Jacobian of a curve.
These questions are trivial for $g=1$.
For $g=2$, the principally polarized question is answered by Howe \cite[Theorem~1.3]{Howe95} for the ordinary case and by Howe-Maisner-Nart-Ritzenthaler \cite{HoweMaisnerNartRitzenthaler08} in general.
The Jacobian question for $g=2$ is resolved by Howe-Nart-Ritzenthaler \cite{HoweNartRitzenthaler09}.
Moreover, Howe \cite[Theorem~1.2]{Howe95} shows that every simple ordinary odd-dimensional isogeny class is principally polarizable. 
The Jacobian question for the supersingular case in dimension $g=3$ is treated by Nart-Ritzenthaler \cite{NartRitzenthaler08}.

\subsection*{Acknowledgments}
This article is a revised version of Jun Jie Lin's master thesis \cite{JunJieLin} which was written at Utrecht University under the supervision of the author. 
Some mathematical corrections have been made.
At that time, the author was supported by NWO grant VI.Veni.202.107 and is currently supported by Marie-Sk{\l}odowska-Curie Actions - Postdoctoral Fellowships 2023 (project 101149209 - AbVarFq).
The author thanks Gunther Cornellissen for encouraging him to write this paper, Jonas Bergstr\"om for comments on a preliminary version of the paper, Boris Shapiro and Vladimir Kostov for suggesting improvements to the exposition in Section~\ref{sec:poly_real_pos_roots}, and the anonymous reviewer for useful suggestions.

\section{Polynomials with only real positive roots}\label{sec:poly_real_pos_roots}\
This section is devoted to the proof Proposition~\ref{prop:real_pos_only}.
This result characterizes when a real polynomial $f(x)$ has only real roots, that is, is hyberbolic, in terms of the roots of its first derivative $f'(x)$.
We expect that Proposition~\ref{prop:real_pos_only} is well-known, but, since we couldn't locate it in the literature, we include a short proof.
A longer but more direct proof is presented in Appendix~\ref{sec:appendix}.
Proposition~\ref{prop:real_pos_only} should be compared with \cite[Lemma]{Smyth84}, where the author assumes also that the roots of $f'(x)$ are distinct.
The author refers to the result as Robinson's method.

We start by recalling two simple properties, whose proofs we omit.
Let $f(x)$ be a complex polynomial and $\gamma$ a complex number.
We denote the order of vanishing of $f(x)$ at $\gamma$ by $\ord_{\gamma}(f(x))$.
If $\ord_\gamma(f(x))=1$ then we say that $\gamma$ is a simple root of $f(x)$.

\begin{lemma}\label{lem:mult_roots_ord}
    Let $f(x)$ be a complex polynomial and $\gamma$ be a complex root of $f(x)$.
    Then $\ord_\gamma(f'(x)) = \ord_\gamma(f(x))-1$.
\end{lemma}

\begin{lemma}\label{lem:mult}
    Let $f(x)$ be a real polynomial of degree $K>1$ with positive leading coefficient.
    Let $\beta$ be a real root of $f'(x)$.
    Then
    \begin{itemize}
        \item if $\ord_\beta(f'(x))$ is even then $\beta$ is an inflection point for $f(x)$;
        \item if $\ord_\beta(f'(x))$ is odd then $\beta$ is either a local minimum or a local maximum for $f(x)$.
    \end{itemize}
\end{lemma}

Given a real polynomial with only real roots, we will use the following description of the set of roots of its first derivative.
\begin{lemma}\label{lemma:real_roots_explicit_descr}
    Let $f(x)$ be a real polynomial of degree $K>1$ with only real roots.
    Denote by and sort as 
    \[ \gamma_s < \ldots < \gamma_1 \]
    the distinct roots of $f(x)$.
    Then $f'(x)$ has only real roots, which are the elements of the set
    \[ \mathfrak{R}\coloneqq\{ \gamma_i : 1\leq i \leq s \text{ and }\ord_{\gamma_i}(f(x))>1 \} \cup
       \{ \delta_{s-1},\ldots,\delta_1 \},
    \]
    where $\delta_{s-1},\ldots,\delta_1$ are distinct simple real roots of $f'(x)$, each one satisfying $\gamma_{i+1}<\delta_i<\gamma_i$.
\end{lemma}
\begin{proof}
    The fact that $\mathfrak{R}$ consists of roots of $f'(x)$ follows from 
    Lemma \ref{lem:mult_roots_ord} and Rolle's Theorem.
    Let $\mathfrak{R}'$ be the multiset consisting of the elements $r \in \mathfrak{R}$ each one counted $\ord_r(f'(x))$-many times.
    Then, using Lemma~\ref{lem:mult_roots_ord}, we get 
    \begin{align*}
        \abs{\mathfrak{R}'} & = \sum_{i=1}^s \ord_{\gamma_i}(f'(x)) + \sum_{i=1}^{s-1}\ord_{\delta_i}(f'(x)) \\
                            & = \sum_{i=1}^s (\ord_{\gamma_i}(f(x))) - s + \sum_{i=1}^{s-1}\ord_{\delta_i}(f'(x))\\ 
                            & \geq \sum_{i=1}^s (\ord_{\gamma_i}(f(x))) - s + s - 1 \\
                            & = \deg(f(x))-1 = \deg(f'(x)).
    \end{align*}
    This show that $\mathfrak{R}$ contains all the roots of $f'(x)$ and that each $\delta_i$ is a simple root of $f'(x)$, thus concluding the proof.
\end{proof}

\begin{proposition}\label{prop:real_pos_only}
    Let $f(x)$ be a real polynomial of degree $K>1$ with positive leading coefficient.
    Assume that all the roots of $f'(x)$ are real.
    Denote them, possibly with repetitions, by
    \[ \beta_{K-1} \leq \ldots \leq \beta_1. \]
    Consider the following condition:
    \begin{equation}\label{min_max}
        \tag{$\diamond$}
        \begin{cases}
            f(\beta_i)\leq 0 & \text{for every $1\leq i \leq K-1$ with $i$ odd},\\
            f(\beta_i)\geq 0 & \text{for every $1\leq i \leq K-1$ with $i$ even}.
        \end{cases}
    \end{equation}
    Then:
    \begin{enumerate}[(A)]
        \item\label{it:A} all complex roots of $f(x)$ are real if and only if \eqref{min_max} holds;
        \item\label{it:B} all complex roots of $f(x)$ are real and non-negative if and only if \eqref{min_max} holds, $\beta_{K-1}\geq 0$ and $(-1)^Kf(0) \geq 0$;
        \item\label{it:C} all complex roots of $f(x)$ are real and positive if and only if \eqref{min_max} holds, $\beta_{K-1}>0$ and $(-1)^Kf(0)>0$.
    \end{enumerate}
\end{proposition}

\begin{proof}
    Consider the statement
    \begin{equation}\label{mult}
        \tag{$\MULT$}
        \text{Every multiple root $\beta$ of $f'(x)$ is a root of $f(x)$ as well.}
    \end{equation}
    Note that if \eqref{mult} holds then $\ord_{\beta}(f(x))=\ord_{\beta}(f'(x))+1>1$.
    If $f(x)$ has only real roots then \eqref{mult} holds by Lemma~\ref{lemma:real_roots_explicit_descr}.
    If we instead assume that \eqref{min_max} holds then \eqref{mult} holds as well: 
    if $\beta_i$ is a multiple root of $f'(x)$, that is, $\beta_{i}=\beta_{i+1}$ or $\beta_{i}=\beta_{i-1}$, then $f(\beta_i)=0$.
    So, we can assume for the rest of the proof that \eqref{mult} holds.

    Now, we start with the proof of \ref{it:A}.
    Because of \eqref{mult}, there is a small enough perturbation of the coefficients of $f(x)$ giving a polynomial $\tilde f(x)$ such the number real roots (counted with multiplicity) of $f(x)$ and $\tilde f(x)$ are the same, the roots $\tilde \beta_i$ of the derivative $\tilde f'(x)$ are distinct, and $\tilde f(\tilde\beta_i)\neq 0$.
    The equivalence in~\ref{it:A} for the polynomial $\tilde f(x)$ is an immediate consequence of Rolle's and the Intermediate Value theorems.
    Since the perturbation described sends a polynomial with only real roots to a polynomial with only real roots, and it is compatible with the inequalities in \eqref{min_max} (it makes them strict), we obtain that~\ref{it:A} holds also for the original polynomial $f(x)$.

    We now prove~\ref{it:B} and~\ref{it:C}.
    By~\ref{it:A}, we can assume that $f(x)$ has only real roots.
    It follows that the two polynomials $f(x)$ are $f'(x)$ are (weakly) interlacing, that is,
    \begin{equation}\label{eq:interlacing}
        \gamma_K \leq \beta_{K-1} \leq \gamma_{K-1} \leq  \cdots \leq \gamma_2 \leq \beta_1\leq \gamma_1.
    \end{equation}
    Then~\ref{it:B} and~\ref{it:C} are easy consequences of the fact that $(-1)^K$ has the same sign of the limit of $f(x)$ at $-\infty$ and Equation~\eqref{eq:interlacing}.
\end{proof}

\section{Polynomials of low degree}\label{sec:poly_low_degree}
In this section, we apply Proposition~\ref{prop:real_pos_only} to polynomials of low degree.
\begin{proposition}\label{prop:deg2}
    Let $f(x) \coloneqq a_2 x^2 + a_1 x + a_0$ be a real polynomial with $a_2 > 0$. Then:
    \begin{itemize}
        \item All complex roots of $f(x)$ are real and non-negative if and only if 
        \[a_1 \geq 0\qquad\text{and}\qquad 0 \leq a_0 \leq \frac{a_1^2}{4 a_2}; \]
        \item All complex roots of $f(x)$ are real and positive if and only if 
        \[a_1 > 0\qquad\text{and}\qquad 0 < a_0 \leq \frac{a_1^2}{4 a_2}.\]
    \end{itemize}
\end{proposition}
\begin{proof}
    The only root of $f'(x)$ is $\beta = - \frac{a_1}{2 a_2} \in \R$.
    The statements are direct consequences of Proposition~\ref{prop:real_pos_only}.
\end{proof}

\begin{proposition}\label{prop:deg3}
    Let $f(x) \coloneqq  a_3 x^3 + a_2 x^2 + a_1 x + a_0$ be a real polynomial with $a_3 > 0$.
    Assume that $a_2^2 - 3 a_1 a_3\geq 0$ and consider the condition:
    \begin{equation}\label{eq:prop:deg3}\tag{$\dagger$}
        \frac{-2 a_2^3 + 9 a_1 a_2 a_3 - 2 (a_2^2 - 3 a_1 a_3)^{3/2} }{27 a_3^2} \leq a_0 \leq \frac{-2 a_2^3 + 9 a_1 a_2 a_3 + 2 (a_2^2 - 3 a_1 a_3)^{3/2} }{27 a_3^2}.
    \end{equation}
    Then:
    \begin{itemize}
        \item All complex roots of $f(x)$ are real and non-negative if and only if 
        $a_2 \geq 0$, $0\leq a_1 \leq \frac{a_2^2}{3a_3}$,
        $a_0 \leq 0$ and \eqref{eq:prop:deg3} holds.
        \item All complex roots of $f(x)$ are real and positive if and only if
        $a_2 > 0$, $0 < a_1 \leq \frac{a_2^2}{3a_3}$,
        $a_0 < 0$ and \eqref{eq:prop:deg3} holds.
    \end{itemize}
\end{proposition}
\begin{proof}
    The derivative $f'(x) = 3 a_3 x^2 + 2 a_2 x + a_1$ has roots
    \[
        \beta_{1} \coloneqq  \frac{- a_2 + \sqrt{ a_2^2 - 3 a_1 a_3}}{3 a_3} \qquad\text{and}\qquad
        \beta_2 \coloneqq  \frac{- a_2 - \sqrt{ a_2^2 - 3 a_1 a_3}}{3 a_3},
    \]
    which are real if and only if $a_2^2 - 3 a_1 a_3\geq 0$. 
    If $f(x)$ has only real roots then the same holds for $f'(x)$ by Lemma~\ref{lemma:real_roots_explicit_descr}.
    Hence, we can assume for the rest of the proof that $a_2^2 - 3 a_1 a_3\geq 0$.
    It follows that $\beta_1 \geq \beta_2$ and 
    \begin{align*}
        f(\beta_1) - a_0 = a_3 \beta_{1}^3 + a_2 \beta_1^2 + a_1 \beta_1 &= \frac{2 a_2^3 - 9 a_1 a_2 a_3 - 2 (a_2^2 - 3 a_1 a_3)^{3/2} }{27 a_3^2}, \\
        f(\beta_2)  - a_0 = a_3 \beta_{2}^3 + a_2 \beta_2^2 + a_1 \beta_2 &= \frac{2 a_2^3 - 9 a_1 a_2 a_3 + 2 (a_2^2 - 3 a_1 a_3)^{3/2} }{27 a_3^2}.
    \end{align*}
    The statements now follow from Propositions~\ref{prop:real_pos_only} and~\ref{prop:deg2} applied to $f'(x)$.
\end{proof}

\begin{proposition}\label{prop:deg4}
    Let $f(x) \coloneqq  a_4 x^4 + a_3 x^3 + a_2 x^2 + a_1 x + a_0$ be a real polynomial with $a_4 > 0$.
    Consider the following inequalities:
    \begin{enumerate}[(i)]
        \item \label{eq:prop:deg4:cond1} $a_3 \leq 0$,
        \item \label{eq:prop:deg4:cond2_star} $0 \leq a_2$,
        \item \label{eq:prop:deg4:cond2} $a_2 \leq \dfrac{3 a_3 ^2}{8a_4},$
        \item \label{eq:prop:deg4:cond3} $a_1 \leq 0$,
        \item \label{eq:prop:deg4:cond4} 
        $ -\dfrac{a_3^3}{8a_4^2} + \dfrac{a_3 a_2}{2a_4} - \dfrac{1}{a_4^2} \left(\dfrac{a_3^2}{4} - \dfrac{2a_2a_4}{3}\right)^{3/2} \leq 
        a_1 \leq 
        -\dfrac{a_3^3}{8a_4^2} + \dfrac{a_3 a_2}{2a_4} + \dfrac{1}{a_4^2} \left(\dfrac{a_3^2}{4} - \dfrac{2a_2a_4}{3}\right)^{3/2}$,
        \item \label{eq:prop:deg4:cond5} $a_0 \leq 0$.
    \end{enumerate}
    Define the real numbers
    \[ 
        u_2 \coloneqq \frac{8 a_2 a_4 - 3 a_3^2 }{16 a_4^2}, \qquad
        u_3 \coloneqq \frac{a_3^3 - 4 a_2 a_3 a_4 + 8 a_1 a_4^2}{32 a_4^3} \qquad\text{and}\qquad
        \Delta \coloneqq u_3^2 + \frac{4}{27} u_2^3,
    \]
    and the set
    \[ S \coloneqq \left\{-\zeta^{k} \omega \left( \frac{9 u_3}{2} + \frac{3\sqrt{\Delta}}{2}\right) - \zeta^{-k} \overline{\omega} \left( \frac{9 u_3}{2} - \frac{3\sqrt{\Delta}}{2} \right)  + \frac{2 u_2^2}{3} : 0\leq k \leq 2 \right\}, \]
    where $\zeta$ is a primitive third root of unity and $\omega$ is a third root of $ \frac{1}{2}(- u_3 + \sqrt{\Delta})$.
    If all complex roots of $f(x)$ are real or \ref{eq:prop:deg4:cond1}--\ref{eq:prop:deg4:cond4} hold then the set $S$ consists of real numbers.
    Assume this is the case, sort the elements of $S$ as $\theta_1 \leq \theta_2 \leq \theta_3$ and consider the condition:
    \begin{equation}\label{eq:prop:deg4}\tag{$\dagger\dagger$}
        \frac{3 a_3^4}{256 a_4^3} - \frac{a_3^2 a_2}{16 a_4^2} + \frac{a_3 a_1}{4 a_4} + a_4 \theta_1  \leq a_0 \leq  \frac{3 a_3^4}{256 a_4^3} - \frac{a_3^2 a_2}{16 a_4^2} + \frac{a_3 a_1}{4 a_4} + a_4 \theta_2.
    \end{equation}
    Then:
    \begin{itemize}
        \item All complex roots of $f(x)$ are real and non-negative if and only if \ref{eq:prop:deg4:cond1}--\ref{eq:prop:deg4:cond5} and \eqref{eq:prop:deg4} hold.
        \item All complex roots of $f(x)$ are real and positive if and only if 
        \ref{eq:prop:deg4:cond1}--\ref{eq:prop:deg4:cond5} and \eqref{eq:prop:deg4} hold, and all inequalities in 
        \ref{eq:prop:deg4:cond1}, \ref{eq:prop:deg4:cond2_star}, \ref{eq:prop:deg4:cond3} and \ref{eq:prop:deg4:cond5} are strict. 
    \end{itemize}    
\end{proposition}
\begin{remark}
    If one assumes \ref{eq:prop:deg4:cond2}, as we do in the statement, the quantity between parentheses in \ref{eq:prop:deg4:cond4} is a non-negative real number.
\end{remark}
\begin{proof}
    By applying Proposition~\ref{prop:deg3} to the first derivative $f'(x) = 4a_4x^3 + 3a_3 x^2 + 2 a_2 x + a_1$,
    we see that all complex roots of $f'(x)$ are real and non-negative if and only if \ref{eq:prop:deg4:cond1}--\ref{eq:prop:deg4:cond4} hold, and that all complex roots of $f'(x)$ are real and positive \ref{eq:prop:deg4:cond1}--\ref{eq:prop:deg4:cond4} hold and all inequalities in 
        \ref{eq:prop:deg4:cond1}, \ref{eq:prop:deg4:cond2_star}, \ref{eq:prop:deg4:cond3} and \ref{eq:prop:deg4:cond5} are strict.

    If $f(x)$ has only real roots then the same holds for $f'(x)$ by Lemma~\ref{lemma:real_roots_explicit_descr}.
    Hence, we can assume for the rest of the proof that the roots of $f'(x)$ are all real.

    We write the roots of $f'(x)$ in terms the coefficients of $f(x)$ using Cardano's method. 
    Consider the depressed cubic 
    \begin{align*}
        f_1(y) \coloneqq \frac{1}{4 a_4} f'\left(y - \frac{a_3}{4 a_4}\right) = y^3 + u_2 y + u_3,
    \end{align*}
    where $u_2$ and $u_3$ are defined as in the statement.
    The real number $\Delta$ defined in the statement is the discriminant of $f_1(y)$.
    Since all the roots of $f'(x)$ are assumed to be real, we have that $\Delta \leq 0$.
    The roots of $f_1(y)$ are $\zeta^k \omega + \zeta^{-k} \overline{\omega}$ for $k=0,1,2$,
    where $\zeta$ is a primitive third root of unity and $\omega \in \C$ is a third root of $ (- u_3 + \sqrt{\Delta})/2$, as in the statement. 
    It follows that the roots of $f'(x)$ are the elements of the set
    \[
        S' \coloneqq \left\{\zeta^k \omega + \zeta^{-k} \overline{\omega} - \frac{a_3}{4 a_4} : 0\leq k \leq 2 \right\} \subseteq \R.
    \]
    Now, consider the set $S$ defined in the statement.
    Since $\Delta \leq 0$, for each $k=0,1,2$, we have
    \begin{equation*}\label{eq:cconj}
        \overline{\zeta^{k} \omega \left( \frac{9 u_3}{2} + \frac{3\sqrt{\Delta}}{2} \right)} = \zeta^{-k} \overline{\omega} \left( \frac{9 u_3}{2} - \frac{3\sqrt{\Delta}}{2} \right).
    \end{equation*}
    This implies that $S$ consists of real numbers.    
    A computation shows that for each $\beta \in S'$ there exists a unique $\theta \in S$ such that
    \[ f(\beta) = a_0 - G - a_4 \theta, \]
    where
    \[ G \coloneqq \frac{3 a_3^4}{256 a_4^3} - \frac{a_3^2 a_2}{16 a_4^2} + \frac{a_3 a_1}{4 a_4}. \]
    Sort the elements of $S'$ as $\beta_3 \leq \beta_2 \leq \beta_1$ and the elements of $S$ as $\theta_1\leq \theta_2 \leq \theta_3$.
    By considering the order of vanishing of $f'(x)$ at each $\beta_i$ and applying Lemma~\ref{lem:mult}, we see that $f(\beta_2) \geq \max\{ f(\beta_1),f(\beta_3)\}$ in all possible cases.
    Since $a_4>0$, we deduce that
    \[  a_0 - G - a_4 \theta_1 = f(\beta_2) 
        \qquad\text{and}\qquad
        a_0 - G - a_4 \theta_2 = \max\{ f(\beta_1),f(\beta_3)\}. \]
    Therefore, \eqref{eq:prop:deg4} is equivalent to $f(\beta_2)\geq 0$, $f(\beta_1)\leq 0$ and $f(\beta_3)\leq 0$.
    The statements are now a direct application of Proposition~\ref{prop:real_pos_only}.
\end{proof}

\begin{remark}\label{rmk:S_unchanged}
    Let notation and hypothesis be as in Proposition~\ref{prop:deg4}.
    Since $\Delta \leq 0$ then
    \[ 
        \overline{\left(\frac{-u_3+\sqrt{\Delta}}{2}\right)\left( \frac{9 u_3}{2} + \frac{3\sqrt{\Delta}}{2}\right)^3} = 
        \left(\frac{u_3+\sqrt{\Delta}}{2}\right)\left( \frac{9\cdot(-u_3)}{2} + \frac{3\sqrt{\Delta}}{2}\right)^3.
    \]
    It follows that the set $S$ is unchanged by replacing $u_3$ with $-u_3$.
\end{remark}

We now apply a reasoning analogous to the one giving Proposition~\ref{prop:deg4} to a polynomial $f(x)$ of degree $5$.
To simplify the exposition, we split the argument in a proposition and a corollary, and assume that $f(x)$ is monic.
This case will suffice to prove Main Theorem~\ref{mainthm:dim5}.
\begin{proposition}\label{prop:deg5monic}
    Let $f(x) \coloneqq  x^5 + a_4 x^4 + a_3 x^3 + a_2 x^2 + a_1 x + a_0$ be a real polynomial.
    Define
    \begin{align*}
        u_2 &\coloneqq \frac{3}{10} a_3 - \frac{3}{25} a_4^2,\\
        u_3 &\coloneqq \frac{2}{125} a_4^3 - \frac{3}{50} a_3 a_4 + \frac{1}{10} a_2, \\
        u_4 &\coloneqq -\frac{3}{625} a_4^4 + \frac{3}{125} a_4^2 a_3 - \frac{2}{25} a_4 a_2 + \frac{1}{5} a_1.
    \end{align*}
    If $u_3 = 0$, let $x_{i_1, i_2}$ for $i_1, i_2 \in \{+1, -1 \}$ be
    \begin{equation}\label{eq:x12_u_3_eq_0}
        x_{i_1, i_2} \coloneqq i_1 \sqrt{- u_2 + i_2 \sqrt{u_2^2 - u_4}}.
    \end{equation}
    If $u_3 \neq 0$, let $x_{i_1, i_2}$ for $i_1, i_2 \in \{+1, -1 \}$ be
    \begin{equation}\label{eq:x12_u_3_neq_0}
        x_{i_1, i_2} \coloneqq \frac{i_1 \sqrt{2y} + i_2 \sqrt{- 4 u_2 - 2y - i_1 \frac{8 u_3}{\sqrt{2y}}}}{2},
    \end{equation}
    where:
    \begin{align*}
        & v_2 \coloneqq - \frac{u_2^2}{3} - u_4, \\
        & v_3 \coloneqq \frac{2 u_2 u_4}{3} - \frac{2 u_2^3}{27} - 2 u_3^2,\text{ and}\\
        & y \coloneqq C - \frac{v_2}{3C} - \frac{2 u_2}{3},\text{ where}\\
        & \text{$C$ is a third root of }
        \begin{cases}
            -v_3 & \text{if $v_2=0$},\\
            \frac{1}{2}\left(- v_3 + \sqrt{v_3^2 + \frac{4}{27} v_2^3}\right) & \text{if $v_2\neq 0$}.
        \end{cases}
    \end{align*}

    Then $x_{i_1, i_2} - \frac{a_4}{5}$, for $i_1, i_2 \in \{+1, -1 \}$ are the four roots of $f'(x)$.
    Assume that they are real, sort them as $\beta_4 \leq \beta_3 \leq \beta_2 \leq \beta_1$, set
    \begin{align*}
        \lambda_1 & \coloneqq \max\left\lbrace \beta_1^4+a_3\beta_1^3+a_2\beta_1^2+a_1\beta_1, 
        \beta_3^4+a_3\beta_3^3+a_2\beta_3^2+a_1\beta_3 \right\rbrace \\
        \lambda_2 & \coloneqq \min\left\lbrace \beta_2^4+a_3\beta_2^3+a_2\beta_2^2+a_1\beta_2, 
        \beta_4^4+a_3\beta_4^3+a_2\beta_4^2+a_1\beta_4 \right\rbrace,
    \end{align*}
    and consider the condition
    \begin{equation}\label{eq:prop:deg5}\tag{$\dagger\dagger\dagger$}
        - \lambda_2 \leq a_0 \leq - \lambda_1.
    \end{equation}
    Then:
    \begin{itemize}
        \item All complex roots of $f(x)$ are real and positive if and only if all the complex roots of $f'(x)$ are real and non-negative, $a_0 \leq 0$ and \eqref{eq:prop:deg5} holds.
        \item All complex roots of $f(x)$ are real and positive if and only if all the complex roots of $f'(x)$ are real and positive, $a_0 < 0$ and \eqref{eq:prop:deg5} holds.
    \end{itemize}
\end{proposition}
\begin{remark}
    The set of the four elements $x_{i_1,i_2}$ for $i_1,i_2\in \{+1,-1\}$ is independent of the choices of third root made in the definition of $C$ and of square roots in Equations~\eqref{eq:x12_u_3_eq_0} and \eqref{eq:x12_u_3_neq_0}, since, as we show in the proof, they are (translates) of the four roots of $f'(x)$.
\end{remark}
\begin{proof}
    In order to determine the bounds on $a_0$, we need to determine the roots of $f'(x) = 5 x^4 + 4 a_4 x^3 + 3 a_3 x^2 + 2 a_2 x + a_1$. 
    We start by making a substitution to obtain a depressed quartic
    \[
        f_1(x) \coloneqq \frac{1}{5} f'\left( x - \frac{a_4}{5} \right) = x^4 + 2 u_2 x^2 + 4 u_3 x + u_4,
    \]
    where $u_2$, $u_3$ and $u_4$ are defined as in the statement.

    If $u_3 = 0$ then $f_1(x)$ is a polynomial in $x^2$.
    Hence, its roots are the elements $x_{i_1,i_2}$ described in \eqref{eq:x12_u_3_eq_0}.

    Now, we consider the case $u_3 \neq 0$.
    For every complex number $y$, the equation
    \begin{equation}\label{eq:Ferrari}
        (x^2 + u_2 + y)^2 = 2 y x^2 - 4 u_3 x + y^2 + 2 y u_2 + u_2^2 - u_4,
    \end{equation}
    has the same set of solutions as the equation $f_1(x) = 0$.
    We choose $y$ in such a way that the right-hand side of \eqref{eq:Ferrari} is a square.
    This happens precisely when the discriminant
    \[ 16 u_3^2 - 4 \cdot 2y \cdot (y^2 + 2 y u_2 + u_2^2 - u_4) \]
    of the quadratic polynomial in $x$ on the right-hand side of \eqref{eq:Ferrari} vanishes, or equivalently, when $y$ is a root of
    \begin{equation}\label{eq:f}
        r(x) \coloneqq x^3 + 2u_2 x^2 + (u_2^2 - u_4) x - 2 u_3^2.
    \end{equation}
    So, assume that $y$ is a real root of $r(x)$.
    Such a $y$ can be explicitly described as
    \[ y = C - \frac{v_2}{3C} - \frac{2 u_2}{3}, \]
    where $v_2$, $v_3$ and $C$ are defined as in the statement.

    Since $u_3\neq 0$, we see that $y\neq 0$ by \eqref{eq:f}.
    Hence, combining the fact that $r(y)/y = 0$ with \eqref{eq:f}, 
    we can rewrite \eqref{eq:Ferrari} as
    \[
        (x^2 + u_2 + y)^2 = \left( x\sqrt{2 y} - \frac{2 u_3}{\sqrt{2y}} \right)^2,
    \]
    which has the same set of solutions of 
    \[
        \left(x^2 + u_2 + y + x\sqrt{2 y} - \frac{2 u_3}{\sqrt{2y}} \right) \left(x^2 + u_2 + y - x\sqrt{2 y} + \frac{2 u_3}{\sqrt{2y}}\right) = 0.
    \]
    So, if $u_3\neq 0$ then roots of $f_1(x)$ are the elements $x_{i_1,i_2}$ described in \eqref{eq:x12_u_3_neq_0}.

    Note that if all the complex roots of $f(x)$ are real, then the same applies to $f'(x)$ by Lemma~\ref{lemma:real_roots_explicit_descr}.
    So, we can assume that this is the case for the rest of the proof.
    Now, observe that \eqref{eq:prop:deg5} is equivalent to $f(\beta_1)\leq 0$, $f(\beta_3)\leq 0$, $f(\beta_2)\geq 0$ and $f(\beta_4)\geq 0$.
    Therefore, the conclusion is an application of Proposition~\ref{prop:real_pos_only}.
\end{proof}

\begin{corollary}\label{cor:monic_deg5}
    Let $f(x) \coloneqq  x^5 + a_4 x^4 + a_3 x^3 + a_2 x^2 + a_1 x + a_0$ be a real polynomial.
    Consider the following inequalities where
    the definitions of $\theta_1$ and $\theta_2$ in \ref{eq:cor:deg5:5}, and of $\lambda_1$ and $\lambda_2$ in \ref{eq:cor:deg5:8} are given below:
    \begin{enumerate}[(i)]
        \item \label{eq:cor:deg5:1} $a_4 \leq 0$,
        \item \label{eq:cor:deg5:2_star} $0 \leq a_3$,
        \item \label{eq:cor:deg5:2} $a_3 \leq \dfrac{2}{5} a_4^2$,
        \item \label{eq:cor:deg5:3}$\dfrac{3}{5}a_3a_4-\dfrac{4}{25}a_4^3 -\dfrac{5}{2}\left( \dfrac{4}{25}a_4^2-\dfrac{2}{5}a_3 \right)^{3/2} \leq a_2 \leq \dfrac{3}{5}a_3a_4-\dfrac{4}{25}a_4^3 + \dfrac{5}{2}\left( \dfrac{4}{25}a_4^2-\dfrac{2}{5}a_3 \right)^{3/2}$,
        \item \label{eq:cor:deg5:4} $a_2 \leq 0$,
        \item \label{eq:cor:deg5:5} $\dfrac{3}{125}a_4^4-\dfrac{3}{25}a_4^2a_3+\dfrac{2}{5}a_2a_4+5\theta_1 \leq a_1 \leq \dfrac{3}{125}a_4^4-\dfrac{3}{25}a_4^2a_3+\dfrac{2}{5}a_2a_4+5\theta_2$,
        \item \label{eq:cor:deg5:6} $0 \leq a_1$,
        \item \label{eq:cor:deg5:7} $a_0 \leq 0$,
        \item \label{eq:cor:deg5:8} $-\lambda_2 \leq a_0 \leq -\lambda_1$.
    \end{enumerate}
    Define $u_2$ and $u_3$ as in Proposition~\ref{prop:deg5monic}.
    Set
    \[ S \coloneqq \left\{-\zeta^{k} \omega \left( \frac{9 u_3}{2} + \frac{3\sqrt{\Delta}}{2}\right) - \zeta^{-k} \overline{\omega} \left( \frac{9 u_3}{2} - \frac{3\sqrt{\Delta}}{2} \right)  + \frac{2 u_2^2}{3} : 0\leq k \leq 2 \right\}, \]
    where $\Delta \coloneqq u_3^2 + \frac{4}{27} u_2^3$, $\zeta$ is a primitive third root of unity and $\omega$ is a third root of $ \frac{1}{2} (- u_3 + \sqrt{\Delta})$.
    If all the roots of $f(x)$ are real or \ref{eq:cor:deg5:1}--\ref{eq:cor:deg5:4} hold then the set $S$ consists of real numbers $\theta_1 \leq \theta_2 \leq \theta_3$.
    Assume this is the case and that $f(x)$ has real roots or that \ref{eq:cor:deg5:5} and \ref{eq:cor:deg5:6} hold, and define $\lambda_1$ and $\lambda_2$ as in Proposition~\ref{prop:deg5monic}.
    Then:
    \begin{itemize}
        \item All the complex roots of $f(x)$ are real and non-negative if and only if \ref{eq:cor:deg5:1}--\ref{eq:cor:deg5:8} hold.
        \item All the complex roots of $f(x)$ are real and positive if and only if \ref{eq:cor:deg5:1}--\ref{eq:cor:deg5:8} hold, and all inequalities in 
        \ref{eq:cor:deg5:1}, \ref{eq:cor:deg5:2_star}, \ref{eq:cor:deg5:4} \ref{eq:cor:deg5:6} and \ref{eq:cor:deg5:7} are strict.
    \end{itemize}
\end{corollary}
\begin{remark}
    The quantity in parentheses in \ref{eq:cor:deg5:3} is a positive real number if \ref{eq:cor:deg5:2} holds.
    Also, the set $S$ in the statement is independent of the choices made in the definitions of $\zeta$ and $\omega$.
\end{remark}
\begin{proof}
    Observe that if all the complex roots of $f(x)$ are real, then the same applies to $f'(x)$ and $f''(x)$ by Lemma~\ref{lemma:real_roots_explicit_descr}.
    Moreover, the inequalities \ref{eq:cor:deg5:1}--\ref{eq:cor:deg5:4} imply that all the complex roots of $f''(x)$ are real by Proposition~\ref{prop:deg3}.
    Also, \ref{eq:cor:deg5:1}--\ref{eq:cor:deg5:6} imply that all the complex roots of $f'(x)$ are real by Proposition~\ref{prop:deg4}.
    Hence, for the rest of the proof, we can assume that all the roots of $f'(x)$ and $f''(x)$ are real.
    Therefore, the statements follow by combining Propositions~\ref{prop:deg3}, \ref{prop:deg4} and \ref{prop:deg5monic}.
\end{proof}

\section{\texorpdfstring{$q$}{q}-Weil polynomials}\label{sec:qWeil}
Let $\Fq$ be a finite field of positive characteristic $p$ with $q=p^n$ elements.
\begin{definition}\label{def:qWeil}
    A $q$-Weil polynomial is a monic polynomial in $\Z[x]$ of even positive degree $2g$ whose set of complex roots has the form $\{w_1,\bar w_1,\ldots,w_g,\bar w_g\}$ and each $w_i$ has absolute value $\sqrt{q}$.
    We denote the set of $q$-Weil polynomials of degree $2g$ by $\Wq{g}$.
\end{definition}

The following lemma contains easy but important observations about $q$-Weil polynomials.
These will be (often implicitly) used throughout the rest of the text.
\begin{lemma}\label{lem:basics}
    Let $h(x)$ be in $\Wq{g}$. Then the following statements hold:
    \begin{enumerate}[(i)]
        \item \label{lem:basics:realroots} all real roots of $h(x)$, if any, occur with even multiplicities.
        \item \label{lem:basics:symm} $h(x)=(x^{2g}/q^g)h(q/x)$.
        \item \label{lem:basics:coeffs} there are integers $a_1,\ldots,a_g$ such that
        \begin{equation*}
            h(x) = x^{2g} + a_1 x^{2g-1} + \cdots + a_{g-1} x^{g+1} + a_g x^{g} + q a_{g-1} x^{g-1} + \cdots + q^{g-1} a_{1} x + q^g.
        \end{equation*}
    \end{enumerate}
\end{lemma}
\begin{proof}
    Statement \ref{lem:basics:realroots} follows immediately from the definition.
    Let $w$ be a complex root of $h(x)$.
    We have $\bar w=q/w$.
    Then $h(x)$ and $(x^{2g}/q^g)h(q/x)$ are both monic polynomials with the same set of roots.
    Hence, they are equal as in~\ref{lem:basics:symm}.
    Finally, \ref{lem:basics:coeffs} follows from \ref{lem:basics:symm}, by comparing the coefficients.
\end{proof}

\begin{example}[$g=1$]\label{ex:g1}
    The polynomial $h(x) \coloneqq  x^2+ax+q$ in $\Z[x]$ is in $\Wq{1}$ if and only if $\abs{a} \leq 2\sqrt{q}$.
    Moreover, $h(x) \in \Wq{1}$ has real roots if and only if $\abs{a} = 2\sqrt{q}$, which can only happen if $q$ is a square.
\end{example}

\begin{remark}\label{rmk:HTtheory}
    In this remark we explain how to determine whether a $q$-Weil polynomial is the characteristic polynomial of Frobenius of some isogeny class of abelian varieties over $\Fq$ of dimension $g$.
    Recall that $q=p^n$.
    We refer to \cite{WaterhouseMilne71} (especially the discussion on page 59) for details.

    If $A$ is an abelian variety over $\Fq$ of dimension $g$, then 
    we have an $\Fq$-isogeny decomposition
    \[ A \sim B_1^{n_1}\times \ldots \times B_r^{n_r}, \]
    where $B_1,\ldots,B_r$ are simple pairwise non-$\Fq$-isogenous abelian varieties and $n_1,\ldots,n_r$ are positive integers.
    It follows from \cite[Theorem 1.(b)]{Tate66} that
    \[ h_A(x) = h_{B_1}(x)^{n_1}\cdots h_{B_r}(x)^{n_r}. \]
    By \cite[Theorem 2.(e) and page~142]{Tate66}, we have that if $A$ is a simple abelian variety over $\Fq$ then
    \begin{equation}\label{eq:h_simple}
        h_{A}(x) = m(x)^e,
    \end{equation} 
    where $m(x)$ a monic $\Q$-irreducible polynomial in $\Z[x]$ and $e$ is the least common denominator of the rational numbers in the set
    \begin{equation}\label{eq:inv}
        \left\{ \frac{v_p(g_1(0))}{n}, \ldots , \frac{v_p(g_s(0))}{n} \right\},
    \end{equation} 
    where $g_1(x),\ldots,g_s(x)$ are the irreducible factors of $m(x)$ over $\Q_p[x]$, $v_p$ is the valuation of $\Q_p$, and we add $1/2$ to the set \eqref{eq:inv} if $m(x)$ has a real root.
    Note that if $m(x)$ has no real roots then it is a $q$-Weil polynomial.
    If instead $m(x)$ has a real root then $m(x)$ is not a $q$-Weil polynomial. Indeed, if $q$ is a square then $m(x) = (x + \sqrt{q})$ or $m(x) = (x - \sqrt{q})$, while if $q$ is not a square then $m(x) = x^2-q$.
    Moreover, $m(x)^2$ is a $q$-Weil polynomial and in fact a characteristic polynomial since $e=2$.

    The work of Honda \cite{Honda68} provides a converse statement.
    If $w\in \C$ is an algebraic integer with complex absolute value $\sqrt{q}$ and minimal polynomial $m(x)$ over $\Q$ then there exists a simple abelian variety $A$ over $\Fq$ with characteristic polynomial as in \eqref{eq:h_simple}.

    Therefore, for 
    \[ h(x) = m_1(x)^{n_1}\cdots m_r(x)^{n_r} \]
    in $\Wq{g}$ with each $m_i(x)$ irreducible over $\Q$, we have that $h(x)$ is a characteristic polynomial of some abelian variety over $\Fq$ of dimension $g$ if and only if, for each $i=1,\ldots,r$, the exponent $n_i$ is a multiple of an exponent $e$ defined as above.
\end{remark}

\begin{example}\label{ex:real_roots}
    Together with our definition of $q$-Weil polynomial, which is used for example in \cite{Haloui10} and \cite{HalouiSingh12}, there is (at least) another definition in the literature in which a $q$-Weil polynomial is a monic integer polynomial whose complex roots have absolute value $\sqrt{q}$.
    To avoid confusion, we will call polynomials satisfying this second more general convention, which is used for example in \cite{LMFDB_art_21}, generalized $q$-Weil polynomials.
    We will not use this more general notion outside this example.
    
    A generalized $q$-Weil polynomial can have odd degree.
    For example, if $q=p^n$ with $n$ odd then the minimal polynomial $m_1(x)$ of $\sqrt{q}$ over $\Q$, which is $m_1(x)=x^2-q$, is a generalized $q$-Weil polynomial but not a $q$-Weil polynomial.
    Also, if $n$ is even then $m_2(x)=x-\sqrt{q}$ and $m_3(x)=(x-\sqrt{q})(x+\sqrt{q})$ are generalized $q$-Weil polynomials but not $q$-Weil polynomial.

    Nevertheless, under the same assumptions on $n$ as above, $m_1(x)^2$, $m_2(x)^2$ and $m_3(x)^2$ are $q$-Weil polynomials.
    In fact, they are characteristic polynomials of abelian varieties over $\Fq$ of dimensions $2$, $1$ and $2$, respectively; see Remark~\ref{rmk:HTtheory}.
\end{example}

\begin{example}\label{ex:Weil_not_char}
    The isogeny class over $\F_4$ with LMFDB-label \href{https://www.lmfdb.org/Variety/Abelian/Fq/4/4/ai_bk_aei_ka}{4.4.ai\_bk\_aei\_ka} has characteristic polynomial $h(x)=m(x)^2$, where $m(x)$ is the irreducible polynomial
    \[ m(x) = x^{4} - 4 x^{3} + 10 x^{2} - 16 x + 16. \]
    By Remark~\ref{rmk:HTtheory}, the polynomial $m(x)$ is a $4$-Weil polynomial, but it is not the characteristic polynomial of an abelian variety over $\F_4$.
\end{example}

The next well-known proposition allows us to describe $\Wq{g}$, which consists of integer polynomials of degree $2g$, in terms of the roots of real polynomials of degree $g$.
We include a proof for completeness. 
\begin{proposition}\label{prop:plus_minus}
    Let $h(x)$ be a polynomial in $\Z[x]$ of the form
    \begin{equation*}
        h(x) = x^{2g} + a_1 x^{2g-1} + \cdots + a_{g-1} x^{g+1} + a_g x^{g} + q a_{g-1} x^{g-1} + \cdots + q^{g-1} a_{1} x + q^g.
    \end{equation*}
    Then $h(x)$ is in $\Wq{g}$ if and only if there exists $\omega_1, \ldots , \omega_g \in \C$ such that
    \begin{equation}\label{eq:prod_quad}
        h(x) = \prod_{i=1}^{g} (x^2 + \omega_i x + q)
    \end{equation}
    and the sets
    \[
        \mathcal{S}^+ \coloneqq  \{ 2 \sqrt{q} - \omega_1, \ldots, 2 \sqrt{q} - \omega_g \}
        \qquad\text{and}\qquad
        \mathcal{S}^- \coloneqq  \{ 2 \sqrt{q} + \omega_1, \ldots, 2 \sqrt{q} + \omega_g \}
    \]
    are subsets of $\R_{\geq 0}$. 
    Moreover, $h(x)$ has no real roots if and only if $0 \notin \mathcal{S}^+ \cup \mathcal{S}^-$.
\end{proposition}
\begin{proof}
    Assume that $h(x)$ is a $q$-Weil polynomial with complex roots $w_1,\bar w_1,\ldots,w_g,\bar w_g$.
    For $i=1,\ldots,g$, set $\omega_i\coloneqq -w_i-\bar w_i$.
    Then $h(x)$ can be written as in \eqref{eq:prod_quad} and 
    $\mathcal{S}^+ \cup \mathcal{S}^- \subset \R_{\geq 0}$.
    If $h(x)$ has a real root, say, $w_i=\pm \sqrt{q}$, then $\omega_i=\mp 2\sqrt{q}$.
    Hence, $0 \in \mathcal{S^+}$ or $0 \in 0 \in \mathcal{S^-}$.
    This proves one direction.
    Now, we prove the converse statement.
    Since 
    $\mathcal{S}^+ \cup \mathcal{S}^- \subset \R_{\geq 0}$,
    each $\omega_i$ is a real number satisfying $\abs{\omega_i} \leq 2 \sqrt{q}$.
    Hence, for each $i$, the complex roots of $x^2 + \omega_i x + q$ are either a pair of conjugate non-real numbers or a single real root with multiplicity $2$.
    In either case, the roots are of the form  $\alpha_i$ and $\bar{\alpha_i}$ and have complex absolute value $\sqrt{q}$.
    Hence, $h(x)$ is a $q$-Weil polynomial.
    Assume in addition that $0 \in \mathcal{S^+}$.
    Then there exists an index $i$ such that $\omega_i=2\sqrt{q}$.
    Then $x^2 + \omega_i x + q$, and hence also $h(x)$, has a real root.
    Similarly, we see that if $0 \in \mathcal{S^-}$ then $h(x)$ has a real root, concluding the proof.
\end{proof}
We use Proposition~\ref{prop:plus_minus}, together with the fact that there are formulas in terms of coefficients for the roots of a real polynomial up to degree $4$, to compute $\Wq{2}$ here in Proposition~\ref{ex:g2} and later $\Wq{3}$, $\Wq{4}$ and $\Wq{5}$ in the following sections, thus proving Main Theorems~\ref{mainthm:dim3}, \ref{mainthm:dim4} and \ref{mainthm:dim5}.
The description of $\Wq{2}$ is certainly well known.
See Section~\ref{sec:intro} for references.

\begin{proposition}[$g=2$]\label{ex:g2}
    Let $a$ and $b$ be integers.
    Consider the integer polynomial 
    \[ h(x) \coloneqq  x^4 + ax^3 + bx^2 + aqx + q^2. \]
    Then $h(x) \in \Wq{2}$ if and only if the following conditions hold:
    \begin{enumerate}[(a)]
        \item \label{ex:g2:1} $\abs{a} \leq 4\sqrt{q}$, and
        \item \label{ex:g2:2} $-2q+2\sqrt{q}\abs{a}\leq b \leq 2q + \dfrac{1}{4}a^2$.
    \end{enumerate}
    Moreover, $h(x) \in \Wq{2}$ has a real root if and only if $\abs{a} = 4\sqrt{q}$ or $b=-2q+2\sqrt{q}\abs{a}$.
    Explicitly, the polynomials $h(x) \in \Wq{2}$ with a real root are $x^4-2qx^2+q^2$ and, if $q$ is a square, the polynomials of the form $(x\pm\sqrt{q})^2h_0(x)$ for $h_0(x)\in \Wq{1}$.
\end{proposition}
\begin{proof}
    Let $h^+(x)$ (resp.~$h^-(x)$) be the monic polynomial whose set of roots is $\mathcal{S^+}$ (resp.~$\mathcal{S^-}$), as defined in Proposition~\ref{prop:plus_minus}.
    Explicitly, we have $h^+(x) = x^2 + b_1^+x + b_0^+$ and $h^-(x) = x^2 + b_1^-x + b_0^-$ where:
    \begin{align*}
        b_1^+ & = -4\sqrt{q}-a, & b_1^- & = -4\sqrt{q}+a, \\ 
        b_0^+ & = 2q+b+2\sqrt{q}a, & b_0^- & = 2q+b-2\sqrt{q}a.
    \end{align*}
    Proposition~\ref{prop:deg2} implies that all complex roots of $h^+(x)$ and $h^-(x)$ are real and non-negative if and only if conditions \ref{ex:g2:1} and \ref{ex:g2:2} hold, and that, if this is the case, $0$ is a root of $h^+(x)$ or $h^-(x)$ if and only if $\abs{a} = 4\sqrt{q}$ or $-2q+2\sqrt{q}\abs{a} = b$.
    Therefore, the characterization of whether $h(x)$ is in $\Wq{2}$, possibly with a real root, in terms of $a$ and $b$ given in the statement is a direct application of Proposition~\ref{prop:plus_minus}.
    
    The polynomials given in the second part of the statement are all in $\Wq{2}$ and have at least a real root.
    We now show that the list is complete.
    Pick $h(x) \in \Wq{2}$ satisfying $\abs{a} = 4\sqrt{q}$ or $b = -2q+2\sqrt{q}\abs{a}$.
    If $a=0$ then $b=2q$, that is, $h(x) = x^4-2qx^2+q^2$.
    Assume for the rest of the proof that $a\neq 0$.
    If $\abs{a}=4\sqrt{q}$ then $q$ is a square and, since $-a$ is the sum of the complex roots of $h(x)$, we must have $h(x)=(x\pm\sqrt{q})^4$ which is the square of a polynomial in $\Wq{1}$.
    Finally, if $\abs{a} \notin \{0,4\sqrt{q}\}$ then not all roots of $h(x)$ are real and $q$ is again a square since   $b$ must be in $\Z$.
    Hence, we can write $h(x) = (x\pm \sqrt{q})^2h_0(x)$ for $h_0(x)\in \Z[x]$ with non-real roots.
    This implies that $h_0(x)$ is in $\Wq{1}$, thus completing the proof.
\end{proof}

\section{Proof of Main Theorem~\texorpdfstring{\ref{mainthm:dim3}}{A}}\label{sec:deg6}
Let $h(x)\coloneqq x^6 + a_1 x^5 + a_2 x^4 + a_3 x^3 + a_2 q x^2 + a_1 q^2 x + q^3$ be as in Main Theorem~\ref{mainthm:dim3}.
Let $h^+(x)$ (resp.~$h^-(x)$) be the monic polynomial whose set of roots is $\mathcal{S^+}$ (resp.~$\mathcal{S^-}$), as defined in Proposition~\ref{prop:plus_minus}.
Explicitly, we have
\[
    h^{+}(x) = x^3 + b_2^+ x^2 + b_1^+ x + b_0^+ \qquad\text{and}\qquad
    h^{-}(x) = x^3 + b_2^- x^2 + b_1^- x + b_0^-,
\]
where:
\begin{align*}
    \begin{array}{lll}
    b_2^+ = - a_1 -6\sqrt{q} , &\quad&b_2^- =  a_1 -6\sqrt{q}, \\
    b_1^+ = a_2 + 4 \sqrt{q} a_1 + 9q, &&b_1^- = a_2 - 4 \sqrt{q} a_1 + 9q, \\
    b_0^+ = - a_3  - 2 \sqrt{q} a_2 - 2q a_1  -2q \sqrt{q}, &&b_0^+ = a_3 - 2 \sqrt{q} a_2 + 2q a_1 -2q \sqrt{q}.
    \end{array}    
\end{align*}
Now, Proposition~\ref{prop:plus_minus} implies that $h(x)$ is in $\Wq{3}$ if and only if all complex roots of $h^+(x)$ and $h^-(x)$ are real and non-negative, which, by Proposition~\ref{prop:deg3}, is equivalent to conditions \ref{mainthm:dim3:a}--\ref{mainthm:dim3:d} from statement.
The characterization in terms of the coefficients of which $h(x)$ is in $\Wq{3}$ has a real root follows  Proposition~\ref{prop:plus_minus} and Proposition~\ref{prop:deg3} as well.

The polynomials listed in \ref{mainthm:dim3:rr0}--\ref{mainthm:dim3:rr3} are all in $\Wq{3}$ and have real roots.
We are left to show that this list contains all $q$-Weil polynomials of degree $6$ with real roots.
Let $h(x)$ be such a polynomial.

Assume first that $q$ is not a square.
Then $\sqrt{q}$ and $-\sqrt{q}$ have minimal polynomial over $\Q$ equal to $(x^2-q)$.
Since, all real roots have even multiplicity by Lemma~\ref{lem:basics}.\ref{lem:basics:realroots}, we deduce that $(x^2-q)^2$ divides $h(x)$.
Hence, we are in case \ref{mainthm:dim3:rr0}.

Assume for the rest of the proof that $q$ is a square.
Again, by Lemma~\ref{lem:basics}.\ref{lem:basics:realroots}, $h(x)$ must have $6$, $4$ or $2$ real roots.
These three possibilities correspond to cases \ref{mainthm:dim3:rr1}, \ref{mainthm:dim3:rr2} and \ref{mainthm:dim3:rr3}, respectively.


\section{Proof of Main Theorem~\texorpdfstring{\ref{mainthm:dim4}}{B}}\label{sec:deg8}
Let $h(x)\coloneqq x^8 + a_1 x^7 + a_2 x^6 + a_3 x^5 + a_4 x^4 + a_3 q x^3 + a_2 q^2 x^2 + a_1 q^3 x + q^4$ be as in Main Theorem~\ref{mainthm:dim4}.
Let $h^+(x)$ (resp.~$h^-(x)$) be the monic polynomial whose set of roots is $\mathcal{S^+}$ (resp.~$\mathcal{S^-}$), as defined in Proposition~\ref{prop:plus_minus}.
Explicitly, we have
\begin{align*}
    h^+(x) & = x^4 + b_3^+ x^3 + b_2^+ x^2 + b_1^+ x  + b_0^+,\\
    h^-(x) & = x^4 + b_3^- x^3 + b_2^- x^2 + b_1^- x  + b_0^-
\end{align*}
where:
\begin{align*}
\begin{array}{lll}
    b_3^+ = - a_1 - 8 \sqrt{q}, &&b_3^- = a_1 - 8 \sqrt{q}, \\
    b_2^+ = a_2 + 6 \sqrt{q} a_1 + 20 q, &&b_2^- = a_2 - 6 \sqrt{q} a_1 + 20 q, \\
    b_1^+ = -a_3 - 4 \sqrt{q} a_2 - 9q a_1 - 16 q \sqrt{q}, &&b_1^- = a_3 - 4 \sqrt{q} a_2 + 9q a_1 - 16 q \sqrt{q}, \\
    b_0^+ = a_4 + 2 \sqrt{q} a_3 + 2q a_2 + 2 q \sqrt{q} a_1 + 2 q^2, &&b_0^- = a_4 - 2 \sqrt{q} a_3 + 2q a_2 - 2 q \sqrt{q} a_1 + 2 q^2.
\end{array}    
\end{align*}
Proposition~\ref{prop:plus_minus} implies that $h(x)$ is in $\Wq{4}$ if and only if all complex roots of $h^+(x)$ and $h^-(x)$ are real and non-negative.
By applying Proposition~\ref{prop:deg4} to both $h^+(x)$ and $h^-(x)$ and manipulating the inequalities, we immediately get \ref{mainthm:dim4:a}--\ref{mainthm:dim4:e} of Main Theorem~\ref{mainthm:dim4}.
Observe that substituting the coefficients of both $h^+(x)$ and $h^-(x)$ in Proposition~\ref{prop:deg4} gives the same set $S$ as discussed in Remark~\ref{rmk:S_unchanged}.
Hence, by substituting the coefficients of $h^+(x)$ or $h^-(x)$, we see that \ref{mainthm:dim4:f} is equivalent to \eqref{eq:prop:deg4} from Proposition~\ref{prop:deg4}.

The characterization in terms of the coefficients of which $h(x)$ is in $\Wq{4}$ has a real root follows Proposition~\ref{prop:plus_minus} and Proposition~\ref{prop:deg4}.

It is clear that the polynomials described in cases \ref{mainthm:dim4:rr0} and \ref{mainthm:dim4:rr1} belong to $\Wq{4}$ and have real roots.
We will show that if $h(x)$ in $\Wq{4}$ with a real root then it is of the form described in case \ref{mainthm:dim4:rr0} or case \ref{mainthm:dim4:rr1}.
If $q$ is a square then, since real roots have even multiplicity by Lemma~\ref{lem:basics}.\ref{lem:basics:realroots}, we see immediately that we are in case \ref{mainthm:dim4:rr0}.
Assume now that $q$ is not a square.
Then the minimal polynomial of $\sqrt{q}$ over $\Q$, which is $x^2-q$, must divide $h(x)$.
In fact, again by Lemma~\ref{lem:basics}.\ref{lem:basics:realroots}, we have that $(x^2-q)^2$ divides $h(x)$.
Hence, we are in case \ref{mainthm:dim4:rr1}.


\section{Proof of Main Theorem~\texorpdfstring{\ref{mainthm:dim5}}{C}}\label{sec:deg10}
Let $ h(x) \coloneqq x^{10} + a_1 x^9 + a_2 x^8 + a_3 x^7 + a_4 x^6 + a_5 x^5 + a_4 q x^4 + a_3 q^2 x^3 + a_2 q^3 x^2 + a_1 q^4 x + q^5$ be as in Main Theorem~\ref{mainthm:dim5}.
Let $h^+(x)$ (resp.~$h^-(x)$) be the monic polynomial whose set of roots is $\mathcal{S^+}$ (resp.~$\mathcal{S^-}$), as defined in Proposition~\ref{prop:plus_minus}.
Explicitly, we have 
\begin{align*}
    h^+(x) & = x^5 +b_4^+ x^4 + b_3^+ x^3 + b_2^+ x^2 + b_1^+ x  + b_0^+,\\
    h^-(x) & = x^5 +b_4^- x^4 + b_3^- x^3 + b_2^- x^2 + b_1^- x  + b_0^-,
\end{align*}
where:
\begin{align*}
    b_4^+ &= - a_1 -10  \sqrt{q}, \\
    b_4^- &= a_1 - 10 \sqrt{q}, \\
    b_3^+ &= a_2 + 8 \sqrt{q} a_1 + 35 q,\\
    b_3^- &= a_2 - 8 \sqrt{q} a_1 + 35 q, \\
    b_2^+ &= -a_3 - 6 \sqrt{q} a_2 - 20q a_1 - 50 q \sqrt{q},\\
    b_2^- &= a_3 - 6 \sqrt{q} a_2 + 20q a_1 - 50 q \sqrt{q}, \\
    b_1^+ &= a_4 + 4 \sqrt{q} a_3 + 9q a_2 + 16 q \sqrt{q} a_1 + 25 q^2,\\
    b_1^- &= a_4 - 4 \sqrt{q} a_3 + 9q a_2 - 16 q \sqrt{q} a_1 + 25 q^2, \\
    b_0^+ &= -a_5 - 2 \sqrt{q} a_4 - 2q a_3 - 2 q \sqrt{q} a_2 - 2 q^2 a_1 - 2 q^2 \sqrt{q},\\
    b_0^- &= a_5 - 2 \sqrt{q} a_4 + 2q a_3 - 2 q \sqrt{q} a_2 + 2 q^2 a_1 - 2 q^2 \sqrt{q}.
\end{align*}
Proposition~\ref{prop:plus_minus} implies that $h(x)$ is in $\Wq{5}$ if and only if all complex roots of $h^+(x)$ and $h^-(x)$ are real and non-negative.
So, we want to deduce Main Theorem~\ref{mainthm:dim5} from Corollary~\ref{cor:monic_deg5}.

Let $x_{i_1,i_2}^+$ (resp.~$x_{i_1,i_2}^-$) for $i_1,i_2 = \pm 1$ be the four elements $x_{i_1,i_2}$ defined by using the polynomial $h^+(x)$ (resp.~$h^-(x)$) in Proposition~\ref{prop:deg5monic}.
Define $u_2^+$, $u_3^+$, $u_4^+$, $u_2^-$, $u_3^-$ and $u_4^-$ analogously.
A computation shows that $u_2^+=u_2^-=u_2$ and $u_4^+=u_4^-=u_4$ while $u_3^+=-u_3^-=-u_3$, where $u_2$, $u_3$ and $u_4$ are as in the statement of Main Theorem~\ref{mainthm:dim5}.
In particular, it follows that the four elements $x_{i_1,i_2}$ defined in the statement of Main Theorem~\ref{mainthm:dim5} coincide with the elements $x_{i_1,i_2}^-$.

Now, inequalities \ref{eq:mainthm:dim5:a}--\ref{eq:mainthm:dim5:g} in Main Theorem~\ref{mainthm:dim5} are directly obtained from conditions \ref{eq:cor:deg5:1}--\ref{eq:cor:deg5:7} in Corollary~\ref{cor:monic_deg5} when applied to both $h^+(x)$ and $h^-(x)$, after observing that the set $S$ in Main Theorem~\ref{mainthm:dim5} remains unchanged if we replace $u_3 = u^-_3$ with $-u_3 = u^+_3$ in its definition; see Remark~\ref{rmk:S_unchanged}.
We now show that \ref{eq:mainthm:dim5:h} in Main Theorem~\ref{mainthm:dim5} is equivalent \ref{eq:cor:deg5:8} in Corollary~\ref{cor:monic_deg5} when applied to both $h^+(x)$ and $h^-(x)$.

The discussion above implies that for each choice of $i_1,i_2$ we have $x_{i_1,i_2}^+ = -x_{-i_1,-i_2}^-$.
In all cases we consider, these elements are real numbers.
Sort the four elements $x_{i_1,i_2}^+$ as $\gamma^+_1\leq\gamma^+_2\leq\gamma^+_3\leq\gamma^+_4$.
Similarly, sort the four elements $x_{i_1,i_2} = x_{i_1,i_2}^-$ as $\gamma_1\leq\gamma_2\leq\gamma_3\leq\gamma_4$, as in the statement of Main Theorem~\ref{mainthm:dim5}.
We have
\begin{equation}\label{eq:gammas}
    \gamma_1 = -\gamma^+_4, \quad \gamma_2 = -\gamma^+_3, \quad \gamma_3 = -\gamma^+_2 \quad\text{and}\quad \gamma_4 = -\gamma^+_1 .
\end{equation}
A straightforward computation shows that
\begin{equation}\label{eq:hphm_comp}
    h^-\left(x-\frac{b_4^-}{5}\right) = -h^+\left(-x-\frac{b_4^+}{5}\right) = H(x) + A - a_5,
\end{equation}
where
\[
    H(x)\coloneqq -x^5 -\frac{10}{3}u_2 x^3 -10u_3x^2 - 5u_4x
\]
and
\[
    A\coloneqq - \frac{4a_1^5}{3125} + \frac{a^3\left(a_2 + 15q \right)}{125} - \frac{a_1^2a_3}{25} -\frac{a_1\left( 3 q a_2 -a_4 +5 q^2\right)}{5} + 2qa_3,
\]
as defined in the statement of Main Theorem~\ref{mainthm:dim5}.
Now, condition \ref{eq:cor:deg5:8} in Corollary~\ref{cor:monic_deg5} applied to both $h^+(x)$ and $h^-(x)$ gives the following eight inequalities:
\begin{equation}\label{eq:hphm_ineq}  
    \begin{aligned}
        h^+\left(\gamma^+_1-\frac{b_4^+}{5}\right) \leq 0, & & h^+\left(\gamma^+_3-\frac{b_4^+}{5}\right) \leq 0, & & h^+\left(\gamma^+_2-\frac{b_4^+}{5}\right) \geq 0, & & h^+\left(\gamma^+_4-\frac{b_4^+}{5}\right) \geq 0, \\
        h^-\left(\gamma_1-\frac{b_4^-}{5}\right) \leq 0, & & h^-\left(\gamma_3-\frac{b_4^-}{5}\right) \leq 0, & & h^-\left(\gamma_2-\frac{b_4^-}{5}\right) \geq 0, & & h^-\left(\gamma_4-\frac{b_4^-}{5}\right) \geq 0.
    \end{aligned}
\end{equation}
By combining \eqref{eq:gammas} and \eqref{eq:hphm_comp}, we deduce that the eight inequalities in \eqref{eq:hphm_ineq} are equivalent to the following four inequalities
\[ 
    H(\gamma_1) + A - a_5 \leq 0, \quad  
    H(\gamma_3) + A - a_5 \leq 0, \quad
    H(\gamma_2) + A - a_5 \geq 0, \quad
    H(\gamma_4) + A - a_5 \geq 0,
\]
which are also equivalent to \ref{eq:mainthm:dim5:h} in Main Theorem~\ref{mainthm:dim5}.

To sum up, $h(x)$ is in $\Wq{5}$ if and only if \ref{eq:mainthm:dim5:a}--\ref{eq:mainthm:dim5:h} hold.
Moreover, again by combining Proposition~\ref{prop:plus_minus} with Corollary~\ref{cor:monic_deg5}, we see that a polynomial $h(x)$ in $\Wq{5}$ has a real root if and only if 
any of the inequalities 
in \ref{eq:mainthm:dim5:a}, \ref{eq:mainthm:dim5:b_star}, \ref{eq:mainthm:dim5:d}, \ref{eq:mainthm:dim5:f} or \ref{eq:mainthm:dim5:g}
is an equality.

The explicit list of $q$-Weil polynomials of degree $10$ has a real root is obtained as follows.
Clearly, the polynomials described in cases \ref{mainthm:dim5:rr0} and \ref{mainthm:dim5:rr1} belong to $\Wq{5}$ and have a real root.
To show the converse, we first assume that $q$ is a square.
Then, since real roots have even multiplicity by Lemma~\ref{lem:basics}.\ref{lem:basics:realroots}, we see immediately that we are in case \ref{mainthm:dim5:rr0}.
Assume now that $q$ is not a square.
Then the minimal polynomial of $\sqrt{q}$ over $\Q$, which is $x^2-q$, must divide $h(x)$.
In fact, again by Lemma~\ref{lem:basics}.\ref{lem:basics:realroots}, we have that $(x^2-q)^2$ divides $h(x)$.
Hence, we are in case \ref{mainthm:dim5:rr1}.






\appendix
\section{A second proof of Proposition~\ref{prop:real_pos_only}}
\label{sec:appendix}

In the proof of Proposition~\ref{prop:real_pos_only} given in Section~\ref{sec:poly_real_pos_roots}, we reduce to the simpler case when the roots of the derivative $f'(x)$ are distinct using a perturbation argument.
Here, we provide a second longer but more direct proof, which might be of independent interest.
The idea behind this proof is to carefully keep track of the contributions given by the multiple roots while scanning the graph of $f(x)$.
In Situation~\ref{notn} below, we introduce the notation and collect a series assumption that we will use for the rest of this section.
\begin{situation}\label{notn}\
    Let $f(x)$ be a real polynomial of degree $K>1$ with positive leading coefficient~$L$.
    Assume that $f'(x)$ has only real roots.
    \begin{enumerate}[(i)]
        \item Let $k$ be the number of distinct roots of $f'(x)$ which are not inflection points for $f(x)$.
        If $k>0$ sort them as
        \[ \alpha_k < \ldots < \alpha_1. \]
        Also, set $\alpha_0\coloneqq +\infty$, $f(\alpha_0)\coloneqq 1$, $\alpha_{k+1}\coloneqq -\infty$ and $f(\alpha_{k+1}) \coloneqq (-1)^K$.
        \item We say that $\SIGN$ holds if $f(\alpha_{i})f(\alpha_{i-1})\leq 0$ for $i=1,\ldots,k+1$.
        \item As in Section~\ref{sec:poly_real_pos_roots}, we say that $\MULT$ holds if each root $\beta$ of $f'(x)$ with $\ord_\beta(f'(x))>1$ is also a root of~$f(x)$.
        \item For $x_1<x_2$ belonging to $\R\cup \{-\infty,+\infty\}$, define
        \[ S(x_1,x_2) \coloneqq \sum_{r \in [x_1, x_2) \cap \R} \left( \ord_r(f(x))-\ord_r(f'(x)) \right). \]
    \end{enumerate}
\end{situation}

\begin{lemma}\label{lem:sit}
    Assuming Situation~\ref{notn}, the following statements hold:
    \begin{enumerate}[(I)]
        \item \label{lem:sit:min_max} 
            Assume $k>0$. 
            For $i=1,\ldots,k$, if $i$ is odd then $\alpha_i$ is a local minimum for $f(x)$, while if $i$ is even then $\alpha_i$ is a local maximum for $f(x)$.
        \item \label{lem:sit:kK_parity} 
            $k$ and $K$ have the opposite parity.
        \item \label{lem:sit:monotone} 
            For each $i=1,\ldots,k+1$, the polynomial $f(x)$ is monotone on the interval $(\alpha_i,\alpha_{i-1})$.
        \item \label{lem:sit:g_no_root}
            For $i=1,\ldots,k+1$: 
            if $f(\alpha_i)f(\alpha_{i-1})\geq 0$ then $f(x)$ has no root in $(\alpha_i,\alpha_{i-1})$; 
            if $f(\alpha_i)f(\alpha_{i-1})< 0$ then $f(x)$ has a unique root in $(\alpha_i,\alpha_{i-1})$ which is either simple or an inflection point.
        \item \label{lem:sit:der_no_roots} 
            Assume that $\MULT$ holds.
            For each $i=1,\ldots,k+1$, if  $f(\alpha_i)f(\alpha_{i-1})\geq 0$ then $f'(x)$ has no root in $(\alpha_i,\alpha_{i-1})$.
        \item \label{lem:sit:re_roots_iff_S} 
            The polynomial $f(x)$ has only real roots if and only if $S(-\infty,+\infty) = 1$.
    \end{enumerate}
\end{lemma}
\begin{proof}
    \ref{lem:sit:min_max} follows from the fact that $L>0$.
    We deduce \ref{lem:sit:kK_parity} from the fact that a root $\beta$ of $f'(x)$ is an inflection point for $f(x)$ if and only if $\ord_\beta(f'(x))$ is even.
    \ref{lem:sit:monotone} is a consequence of the fact that on each interval considered there aren't local minima nor local maxima for $f(x)$.
    \ref{lem:sit:g_no_root} is a consequence of \ref{lem:sit:monotone}.
    We prove \ref{lem:sit:der_no_roots} by contradiction: assume that $\beta$ is a root of $f'(x)$ in $(\alpha_i,\alpha_{i-1})$. 
    Then $\beta$ is an inflection point of $f(x)$ and hence, by $\MULT$, is a root of $f(x)$. 
    This contradicts \ref{lem:sit:g_no_root}.
    Finally, since we are assuming that $f'(x)$ has only real roots, we get \ref{lem:sit:re_roots_iff_S}.
\end{proof}

\begin{proposition}\label{prop:real_roots_only}
    Assuming Situation~\ref{notn}.
    Then:
    \begin{enumerate}[(a)]
        \item \label{real_only} The polynomial $f(x)$ has only real roots if and only if $\SIGN$ and $\MULT$ hold.
        \item \label{real_non-neg_only} All complex roots of $f(x)$ are real and non-negative if and only if $\SIGN$ and $\MULT$ hold, all roots of $f'(x)$ are non-negative and $(-1)^K f(0)\geq 0$.
        \item \label{real_positive_only} All complex roots of $f(x)$ are real and positive if and only if $\SIGN$ and $\MULT$ hold, all roots of $f'(x)$ are positive and $(-1)^K f(0)>0$.
    \end{enumerate}
\end{proposition}
\begin{proof}
    We start by proving~\ref{real_only} of the theorem.
    By Lemma~\ref{lemma:real_roots_explicit_descr}, we have that
    \begin{equation}\label{rreq00}
        \text{if $f(x)$ has only real roots then $\MULT$ holds}.
    \end{equation}
    Hence, for the rest of the proof of \ref{real_only}, we can and do assume that $\MULT$ holds.

    Assume first that $k=0$, that is, all roots of $f'(x)$ are inflection points for $f(x)$.
    By Lemma~\ref{lem:sit}.\ref{lem:sit:monotone}, the polynomial $f(x)$ is monotone.
    Hence, there exists a unique $\beta \in \R$ such that $f(\beta)=0$. $\MULT$ implies that $\beta$ is the unique inflection point of $f(x)$. 
    It follows that $f(x)=L(x-\beta)^K$ and $f'(x)=KL(x-\beta)^{K-1}$.
    This concludes the proof of~\ref{real_only} for the case $k=0$.

    Assume now that $k\geq 1$.
    We keep this assumption until the end of the proof of Part~\ref{real_only}.
    Set $S_{k+1}\coloneqq S(-\infty,\alpha_k)-1$ (the $-1$ is added to simplify the exposition below) and, for $i=1,\ldots,k$, set $S_i\coloneqq S(\alpha_i,\alpha_{i-1})$.
    Using this notation we have 
    \[
        S(-\infty,+\infty) = (S_{k+1}+1) + S_{k-1} + \ldots + S_1 .
    \]
    If follows from Lemma~\ref{lem:sit}.\ref{lem:sit:re_roots_iff_S} that 
    \begin{equation}\label{rreq0}
        \sum_{i=1}^{k+1} S_i = 0 \text{ if and only if $f(x)$ has only real roots.}
    \end{equation}
    For each $i=1,\ldots,k+1$, exactly one of the following cases occurs:
    \begin{enumerate}[(i)]
        \item \label{1} If $f(\alpha_i)f(\alpha_{i-1}) > 0$ then $S_i=-1$ : 
            Indeed, by Lemma~\ref{lem:sit}.\ref{lem:sit:g_no_root}, the polynomial $f(x)$ has no root in $(\alpha_i,\alpha_{i-1})$ and, by assumption, $f(\alpha_i)\neq 0$.
            By Lemma~\ref{lem:sit}.\ref{lem:sit:der_no_roots}, also $f'(x)$ has no root on $(\alpha_i,\alpha_{i-1})$.
            So, if $i=k+1$ then we get $S(-\infty,\alpha_k)=0$ and therefore $S_{k+1}=-1$.
            If instead $i\leq k$ then $\alpha_i$ is a simple root of $f'(x)$ by $\MULT$, and hence $S_i = S(\alpha_i,\alpha_{i-1}) = -\ord_{\alpha_i}(f'(x)) = -1$.
        \item \label{2} If $f(\alpha_i)f(\alpha_{i-1}) < 0 $ then $S_i=0$ : 
            By Lemma~\ref{lem:sit}.\ref{lem:sit:g_no_root}, there exists a unique root $\gamma$ of $f(x)$ in $(\alpha_i,\alpha_{i-1})$ which is either simple or an inflection point for $f(x)$.
            In either case, $\ord_{\gamma}(f(x))-\ord_{\gamma}(f'(x)) = 1$.
            So, if $i=k+1$ then we get $S_{k+1}=S(-\infty,\alpha_k)=1-1=0$.
            If $i\leq k$ then, since $\alpha_i$ is a simple root of $f'(x)$ by $\MULT$, we get $S_i=S(\alpha_i,\alpha_{i-1}) = \ord_{\alpha_i}(f(x))-\ord_{\alpha_i}(f'(x)) + \ord_{\gamma}(f(x))-\ord_{\gamma}(f'(x)) = -1 + 1 = 0$.
        \item \label{3} If $f(\alpha_i)f(\alpha_{i-1}) = 0$ then, since we cannot have $f(\alpha_i)=f(\alpha_{i-1})=0$ by Lemma~\ref{lem:sit}.\ref{lem:sit:min_max}, exactly one of the following two possibilities occurs:
            \begin{enumerate}[label=(\roman{enumi}.\alph*), ref=(\roman{enumi}.\alph*)]
                \item \label{3a} $i \leq k$, $f(\alpha_i)=0$ and $f(\alpha_{i-1})\neq 0$. 
                    The only root of $f(x)$ in $[\alpha_i,\alpha_{i-1})$ is $\alpha_i$ by $\MULT$ and Lemma~\ref{lem:sit}.\ref{lem:sit:g_no_root}. 
                    It follows that $S_i=1$ by Lemma~\ref{lem:mult_roots_ord} and Lemma~\ref{lem:sit}.\ref{lem:sit:der_no_roots}.
                    Moreover, $f(\alpha_{i+1})\neq 0$.
                \item \label{3b} $i\geq 2$, $f(\alpha_i) \neq 0$ and $f(\alpha_{i-1}) = 0$.
                    By Lemma~\ref{lem:sit}.\ref{lem:sit:g_no_root} there are no roots of $f(x)$ in $[\alpha_i,\alpha_{i-1})$.
                    Hence, $S_i=-1$ by Lemma~\ref{lem:sit}.\ref{lem:sit:der_no_roots}.
                    Moreover, $f(\alpha_{i-2})\neq 0$.
            \end{enumerate}
            In particular, combining this information, we see that either $(S_{i+1},S_i)=(-1,1)$ (case \ref{3a}), or $(S_i,S_{i-1})=(-1,1)$ (case \ref{3b}).
    \end{enumerate}
    Now, we deduce from \ref{1}--\ref{3} two facts.
    Firstly, we see that $S_{k+1}\in \{-1,0\}$ and that for $i=1,\ldots,k$, we have $S_i \in \{-1,0,1\}$.
    Secondly, for $i=1,\ldots,k$, if $S_i=1$ then $S_{i+1}=-1$.
    Therefore
    \begin{equation}\label{non-pos}
        \sum_{i=1}^{k+1} S_i \leq 0.
    \end{equation}
    By definition, $\SIGN$ does not hold if and only if \ref{1} occurs for some index $i=1,\ldots,k+1$.
    Assume that $\SIGN$ holds.
    For each $i=2,\ldots,k+1$ with $S_i=-1$ we are necessarily in \ref{3b} and hence $S_{i-1}=1$.
    By \eqref{non-pos}, we get $S_{k+1} + \ldots + S_1 = 0$.
    Assume that $\SIGN$ does not hold.
    Let $i$ be an index for which case \ref{1} occurs, so that $S_i=-1$.
    Consider the sequence $(S_{k+1},\ldots,S_1)$.
    Every $1$ occurs only in case \ref{3a}, which always follows a $-1$ from case \ref{3b}.
    Hence, the $-1$ at index $i$ is not balanced by a $1$.
    It follows that $S_{k+1} + \ldots + S_1 < 0$.
    We have shown that
    \begin{equation}\label{notSIGN}
        \sum_{i=1}^{k+1} S_i < 0 \text{ if and only if $\SIGN$ does not hold.}
    \end{equation}
    Therefore, combining \eqref{rreq00}, \eqref{rreq0}, \eqref{non-pos} and \eqref{notSIGN} we obtain that $f(x)$ has only real roots if and only if $\SIGN$ and $\MULT$ hold, concluding the proof of~\ref{real_only} of the theorem.
    \smallskip
    
    We now prove~\ref{real_non-neg_only} and \ref{real_positive_only}.
    We assume that all the roots of $f(x)$ are real, which, as we have shown, is equivalent to have $\SIGN$ and $\MULT$.
    Let $\beta$ be the smallest root of $f'(x)$.
    Hence, $f'(x)$ has only non-negative (resp.~positive) roots if and only if $\beta\geq 0$ (resp.~$\beta>0$).
    By Lemma~\ref{lemma:real_roots_explicit_descr}, there exists a unique root $\gamma$ of $f(x)$ with $\gamma \leq \beta$.
    In particular, $f(x)$ has only non-negative (resp.~positive) roots if and only if $\gamma\geq 0$ (resp.~$\gamma>0$).
    
    Assume that $\gamma \geq 0$.
    By Lemma~\ref{lem:sit}.\ref{lem:sit:monotone}, $f(x)$ is monotone on $(-\infty,\beta)$ and hence on $(-\infty,\gamma)$.
    Hence, if $\gamma>0$ then $(-1)^K$ and $f(0)$ have the same sign, that is, $(-1)^K f(0)>0$. 
    Hence, since $\gamma>0$ trivially implies $\beta>0$, we have proven the `left-to-right' implication of Part~\ref{real_positive_only}.
    Since $\gamma=0$ implies that $\beta\geq 0$ and $(-1)^K f(0)=0$, we obtain also the `left-to-right' implication of Part~\ref{real_non-neg_only}.

    Assume now that $\beta \geq 0$.
    By Lemma~\ref{lem:sit}.\ref{lem:sit:monotone}, $f(x)$ is monotone on $(-\infty,\beta)$ and hence on $(-\infty,0)$.
    So, if $(-1)^K f(0)>0$ then $f(x)$ has no root in $(-\infty,0]$.
    This implies that $\gamma>0$.
    This completes the proof of Part~\ref{real_positive_only}.
    If instead $(-1)^K f(0)=0$ then we must have $\gamma=0$, thus completing the proof of Part~\ref{real_non-neg_only}.
\end{proof}
To directly apply Proposition~\ref{prop:real_roots_only} to a real polynomial $f(x)$ one needs to compute the roots of $f'(x)$, sort them and determine which one are local minima and local maxima.
This last step can be avoided by using Proposition~\ref{prop:real_pos_only}.
\begin{proof}[Second proof of Proposition~\ref{prop:real_pos_only}]
    The assumptions on $f(x)$, $f'(x)$ and $K$ are compatible with Situation~\ref{notn}, from which we borrow the notation.
    By Proposition~\ref{prop:real_roots_only}, in order to prove all statements, it suffices to show that \eqref{min_max} is equivalent to have both $\SIGN$ and $\MULT$.
    By Lemma~\ref{lem:mult}, for each $1\leq i \leq K-1$, if $i$ is odd then $\beta_i$ is either an inflection point or a local minimum for $f(x)$, while if $i$ is even then $\beta_i$ is either an inflection point or a local maximum for $f(x)$.
    In particular, \eqref{min_max} implies $\SIGN$.

    We now show that \eqref{min_max} implies $\MULT$.
    Let $\beta$ be a multiple root of $f'(x)$, so that $K>2$.
    There exists an index $i$ such that $\beta=\beta_i=\beta_{i+1}$.
    Hence, \eqref{min_max} implies that
    \[ f(\beta) \leq \max\{ f(\beta_{i}) : i\ \text{odd}\}\leq 0 \leq \min\{ f(\beta_{i}) : i\ \text{even}\} \leq f(\beta). \]
    Therefore, $f(\beta)=0$, as required.

    Now, on the one hand, if $\beta_i$ is a local minimum or a local maximum and $\SIGN$ holds then $\beta_i$ satisfies the condition stated in \eqref{min_max}.
    On the other hand, if $\beta_i$ is an inflection point and $\MULT$ hold then $\beta_i$ satisfies the condition stated in \eqref{min_max}.
    Therefore, if $\SIGN$ and $\MULT$ hold then \eqref{min_max} holds.
\end{proof}

\bibliographystyle{plain}
\renewcommand{\bibname}{References} 
\bibliography{references} 

\end{document}